\documentclass[12pt]{article}
\usepackage[utf8]{inputenc}
\usepackage{graphicx}
\usepackage{caption}
\usepackage{latexsym}
\usepackage[english]{babel}

\usepackage{amscd,amsfonts,amsmath,amssymb,amsthm}
\usepackage[mathscr]{eucal}

\usepackage[usenames]{color}

\input{xy}
\xyoption{all}

\usepackage{tensor} 
\usepackage{mathdots} 

\usepackage[toc,page]{appendix} 
\usepackage{hyperref}
\hypersetup{
    colorlinks=true,
    linkcolor=blue,
    filecolor=magenta,
    urlcolor=magenta,
    citecolor=blue
}

\newcommand{\hb} {\rule{5in}{0.1pt}}
\newcommand{\R}{\mathbb{R}}

\theoremstyle{plain}
\newtheorem{theorem}{Theorem}

\newtheorem{lemma}{Lemma}
\newtheorem{cor}{Corollary}

\theoremstyle{definition}

\theoremstyle{remark}
\newtheorem{remark}{Remark}

\title{
On the asymptotics and the non-holonomic character \\
of first returns in the standard Euclidean lattice
}
\author{Dorin Dumitra\c{s}cu\footnote{Corresponding author: \texttt{ddumitrascu@adrian.edu}} and Liviu Suciu}
\date{\today}

\begin{document}

\maketitle

\begin{abstract}

We give precise asymptotics to the number of first time returning random walks in the standard orthogonal lattice in $\R^d$ and we prove that these numbers do not form a $P$-recursive sequence.
In the process, the known asymptotics of the number of closed walks are obtained in an elementary way, by using a combinatorial and geometric multiplication principle together with the classical theory of Legendre polynomials.
By showing that the relevant generating functions are $G$-functions, we use a form of the Hadamard convolution to find their singularities in all dimensions and give the ODEs that they satisfy for $d\leq 5$, some of which seem to be new.
We use the Lucas property of the number of closed walks to prove that the corresponding generating function is not invertible as a $G$-function, which immediately implies that the generating function of the first time returning walks is not holonomic.
We propose a few conjectures on the form of the asymptotic coefficients and of the ODEs.

\end{abstract}

\vspace{10pt}

\vspace{7pt}


\section{Introduction}
Random walks have been studied extensively for more than a century in connection to Brownian motion, diffusion, behavior of financial time series, and that of many other stochastic processes.
In this paper, we think of random walks as occurring in an Euclidean space $\mathbb R^d$, along the directions of a fixed orthogonal system of coordinates, with equally likely probability of motion in the direction of the axes.
We present some interesting connections that random walks afford between combinatorics, theory of Legendre polynomials, $G$-functions, and singularity analysis.

The first part of our paper has combinatorial and geometric themes and presents two multiplication principles satisfied by the number of walks ending at various points in $\mathbb R^d$.
Theorem~\ref{T:layers-allD} describes an inductive, on $d$, procedure to calculate the number of walks of length $n$ that end up at an arbitrary point in $\mathbb R^d$.
We denote by $A_{2n}^{(d)}$ the number of walks of length $2n$ in $\mathbb R^d$ that return at the origin at time $t=2n$. These are known in the literature as ``closed walks" or ``excursions." Similarly, we denote by $B_{2n}^{(d)}$ the number of walks of length $2n$ in $\mathbb R^d$ that return {\it for the first time} at the origin at time $t=2n$.
As part of our results in Section~2, we obtain in Corollary~\ref{C:As-all-d} the following recurrence relation:
$$
A_{2n}^{(d+1)} = \binom{2n}{n} \cdot \sum_{k=0}^{n} \binom{n}{k}^2 \,
        \frac{A_{2k}^{(d)}}{\binom{2k}{k}}.
$$

The main results of the paper are contained in Section~3. In Subsection~\ref{subsec:As}, Theorem~\ref{T:TheoremA}, we reobtain the asymptotics
$$
A_{2n}^{(d)}=a_d\frac{(2d)^{2n}}{(\pi n)^{d/2}} \left(1+ \frac{a_1(d)}{n}+\frac{a_2(d)}{n^2}+\dots +\frac{a_m(d)}{n^m}+
O_{d,m}\left( \frac{1}{n^{m+1}} \right) \right),
$$
where $a_d=d^{d/2}/2^{d-1}$, $a_1(d)= -d/8$, $a_2(d)=(2d^3-3d^2+4d)/384$, while $a_3(d), a_4(d)$ are polynomials in $d$ of degree $4$ and $6$, and we conjecture that all $a_m(d)$ are polynomials in $d$ of degree $[\frac{3m}{2}]$ and divisible by $d$.

These results are proved in a relatively elementary way using the theory of Legendre polynomials and standard results about integrals and Taylor series of basic functions like $\exp$ and  $\log$, while using from Section~2 only the above recurrence and the initial conditions for $d=1$.

In Subsection~\ref{subsec:Bs} we use the theory of holonomic functions and the Tauberian theory of singularity analysis \cite{FS} to obtain asymptotics for the first return paths in the form

$$
B_{2n}^{(d)} = b_d\frac{(2d)^{2n}}{(\pi n)^{d/2}} \left(1+ \frac{b_1(d)}{n}+\frac{b_2(d)}{n^2}+\dots +\frac{b_m(d)}{n^m}+
O_{d,m}\left( \frac{1}{n^{m+1}} \right) \right), \text{ for } d \ge 3 \text{ odd}.
$$

\noindent
Here, $b_d=a_d/m_d^2$, where, for all $d \ge 3$, $m_d$ is the expected number of returns to the origin, plus one:
$
m_d = 1+\sum_{n=1}^{\infty} A_{2n}^{(d)}/(2d)^{2n}<\infty.
$
The case $d$ even additionally involves polynomials in $\log n$. See Theorem~\ref{T:TheoremB} for full details. Using the auxiliary sequence $x_n^{(d)}=A_{2n}^{(d)}/{2n \choose n}$ and applying the theory of Hadamard convolution of analytic functions, we show that the generating functions $F_d$ and $A_d$, of $(x_n)_n$ and $(A_{2n})_n$, satisfy Fuchsian ODEs, which we conjecture to be of degree $(d-1)$ and $d$, respectively. We determine the (finitely many) singularities of these generating functions and their type. For example, $A_d$ has finite order singularities at $1/(2k)^2$, for all $1 \le k \le d$ and $k$ having the same parity as $d$, and at $\infty$. However the generating function $B_d$ is not holonomic for $d\geq 2$, but it is $\Delta$-analytic and the theory of singularity analysis still apples. See Theorems~\ref{T:TheoremXs}, \ref{T:TheoremAs}, and \ref{T:TheoremBs}.
We also obtain, in Theorem~\ref{T:TheoremC}, asymptotic expansions for the (normalized) generating functions in terms of the standard functions $T(w)=1-w$ and $L(w)=-\log (1-w)$, $|w|<1$.

In Subsection~\ref{subsec:AnP} we present the necessary proofs, which are based on the identity (see equation~(\ref{eq:ConvExpN}))
$$
F_{d+1}(w) = \frac{1}{2\pi i}\int_{|t|=1/s_d}
    F_d(t)P(t,w)^{-1/2}dt, \text{ for } P(t,w)=w^2-wt(2w+2)+t^2(w-1)^2.
$$
(Here $s_d > d^2$.) In order to adapt to our case the determination of singularities of Hadamard convolution of two functions, we need to carefully analyze
$$
t_1(w)=\frac{w}{(\sqrt w +1)^2} \text{ and }
t_2(w)=\frac{w}{(\sqrt w -1)^2},
$$
which are the zeros of $P(t,w)$.
The proof of above identity is based on properties of the Legendre polynomials in the complex domain, and especially on the relation
$$
\sum_{m=0}^{\infty} {n+m \choose m}^2 w^m=(1-w)^{-n-1}P_n\left( \frac{1+w}{1-w} \right),
\text{ for } |w|<1,
$$
which seems to have a long history, appearing as early as 1303 (in the form of the equality between the coefficients of $w^k$ on both sides) in the works of the Chinese mathematician Shih-Chieh Chu (\cite{Takacs}).

In Section~\ref{sec:ODEs}, we look at the generating functions $F_d$ and $A_d$ in the context of the theory of holonomic power series and in particular the theory of $G$-functions (which are holonomic power series with coefficients algebraic numbers which satisfy some growth conditions, that automatically happen when they are integers). We prove that $F_d$, for $d\geq 3$, $A_d$ and $A'_d/A_d$, for $d\geq 2$, are not algebraic functions, with the later implying that $1/A_d$ and $B_d$ are not holonomic. Hence the sequence $(B_{2n}^{(d)})_n$ is not $P$-recursive, in contrast with the sequences $(x_{n}^{(d)})_n$ and $(A_{2n}^{(d)})_n$. While most of these follow using the asymptotics proved in Section~3 and general considerations, this method would require the (unknown but likely) transcendentality of $m_{2d+1}\pi^d$ to prove the non-algebraicity of $A'_{2d+1}/A_{2d+1}$, and consequently the non-holonomicity of $B_{2d+1}$. However, both the sequences $( x_n^{(d)} )_n$ and $( A_{2n}^{(d)} )_n$ satisfy the Lucas property \cite{McIntosh}, $u_{np+q}=u_nu_q$ mod $p$, for all $n \ge 0$, $0 \le q \le p-1$ and $p$ prime. This allows us to give a direct proof that $1/F_d$, $d \ge 3$, and $B_d$, $d \ge 2$, are not holonomic. We in fact prove that a formal power series $G$ with integral coefficients satisfying Lucas property and for which $G'/G$ is algebraic over $\mathbb Q$, must satisfy $G'/G$ be actually rational over $\mathbb Q$.
The paper ends with an enumeration of the $P$-recurrences and the ODEs satisfied by the sequences $(x_{n}^{(d)})_n$ and $(A_{2n}^{(d)})_n$, for $n=1$ to 5, of which some seem to be new.

This article is a drastically modified version of our preprint \cite{DS}. We thank the anonymous referee for making us aware of some already published results which our preprint rediscovered or which went beyond our initial approach.  We refer here especially to \cite{Joyce73} and to the concept of holonomicity \cite[Appendix B.4]{FS}.

\section{Simple random walks in Euclidean lattices}

To set up the discussion, consider the Euclidean space $\R^d$, with a fixed system of coordinates. We will use the lattice with vertices the points with integer coefficients and with each vertex being joined to the nearest neighbors by an undirected segment.
By a {\it (simple) random walk of length $n$} we understand the motion of a point which starts at the origin $P_0$ and ``walks" along $n$ segments of the lattice, at each vertex choosing with equal probability $1/(2d)$ which segment to follow next.
Sometimes we refer to the number of traveled segments as {\it time}, because this set-up corresponds to a discrete-time and discrete-space walk, one step of the walk occurring with each increment of the time $t$.
In this section we obtain explicit formulas for the {\it number} of random walks of length $n$ that end up at a particular point in $\R^d$, when the starting point is fixed. In other words, we determine the geometric distribution of these numbers. See Figure~\ref{fig:diamonds} below, as an exemplification in $\R^2$. Such formulas have been obtained for some time by Domb \cite{Domb}, using probability distribution functions. Our approach is more geometric and leads to the recurrence formula (\ref{eq:As-all-d}) which allowed us to make connections with Legendre polynomials.

For simplicity, we assume that the walker starts at the origin of the coordinate system. In the case $d=1$, it is easy to see that the number of walks of length $n$ is given by the binomial coefficients. When $d=2$ one obtains a diamond shaped distribution.
Figure~\ref{fig:diamonds} shows this distribution in the case $n=3$ and $n=4$. Each green (darker) dot represents the end of a path. The number next to a dot represents the number of paths that end at that green dot, among all possible walks of length $n$. The gray (lighter) dots are the ends of paths from the previous time step.

\begin{figure}[h]
\begin{center}
    \scalebox{.35}[.4]{\includegraphics{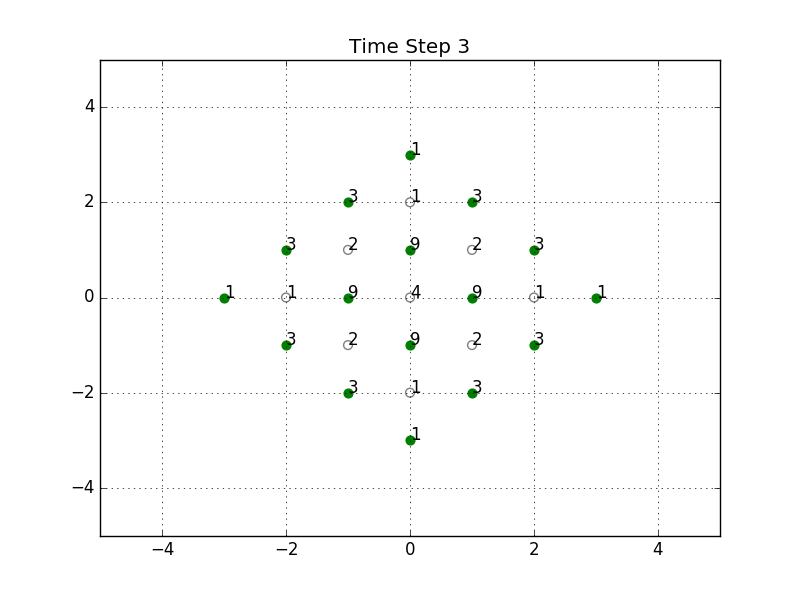}}
    \scalebox{.35}[.4]{\includegraphics{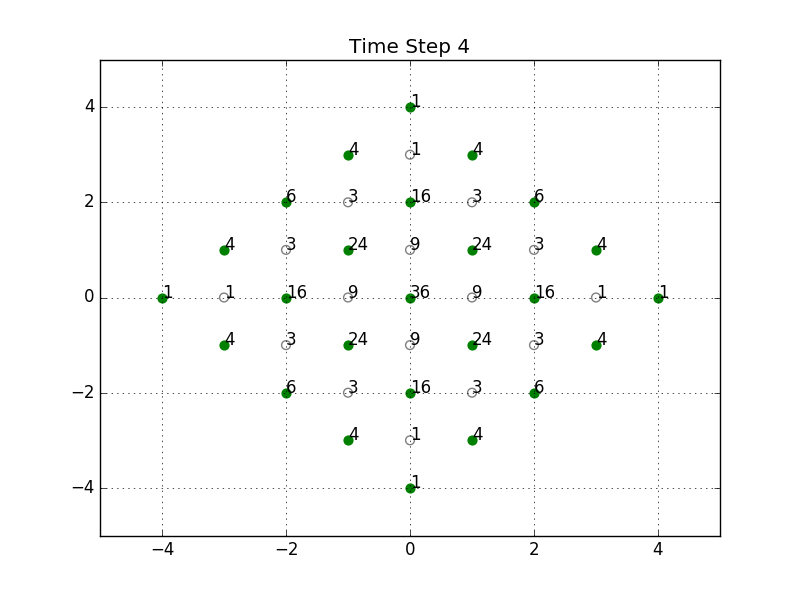}}
\end{center}
 \caption{The diamond shaped distribution of random walks starring at $P_0=(0,0)$.}
 \label{fig:diamonds}
\end{figure}

We refer to this pattern as the {\it first multiplication principle}:

\begin{lemma}\label{L:MP}
The number of random walks of length $n$ in the plane that start at the origin $(0,0)$ and end at point $(k,l)$ is given by
$
\displaystyle{
\binom{n}{\frac{n+k-l}{2}} \binom{n}{\frac{n+k+l}{2}}
}.
$
Here, $-n \le k \le n$, $-n \le l\le n$, $|k| + |l| \le n$, and $n \equiv k+l \pmod{2}$.
\end{lemma}

The proof is by induction on $n$. In its essence, Lemma~\ref{L:MP} is due to the fact that a random walk in $\R^2$ can be viewed as the product of two independent one-dimensional random walks. This is clear if one projects the walk on the orthogonal system consisting of the diagonals of the standard system of coordinates. It is mentioned in the original paper of P\'{o}lya \cite{Polya}, in \cite[Example 12]{FFK}, \cite[Example VI.14]{FS}, and, in the disguised form of generating functions, in \cite[page 323]{Finch}.

%
The three-dimensional case is slightly more complicated. There is no equivalent of Lemma~1, perhaps not surprisingly because the $l^1$ unit ball is not a cube, but a different multiplication principle still holds. The end points of all possible random walks of a fixed length form octahedra and a multiplication principle which holds on the {\it faces} of such octahedra turns out to hold for ``layers" too (see below). Moreover, the result that we were able to first observe for $\R^3$ turns out to hold for any $\R^d$, with $d\geq 3$. See Theorem \ref{T:layers-allD}.

Before analyzing walks, we discuss a sequence of triangles of numbers
$( \mathbf{\tau_n} )_{n\geq 0}$ which might be of interest in its own right.
This sequence is built inductively using the following procedure. At stage $n=0$, $\mathbf{\tau_0}$ consists just of the number 1. At any subsequent stage $n+1$, the numbers from $\mathbf{\tau_n}$ duplicate one slot above, to the left, and to the right, and then they are removed from the triangle. The new triangle $\mathbf{\tau_{n+1}}$ consists of the sums of the numbers that so result at each location. One could interpret $\mathbf{\tau_n}$ as the distributions of the number of random walks of length $n$ in the plane, where the moves happen with equal probability equal to $1/3$ to the left, right, and upward directions only. The diagram bellow shows the first five triangles of this sequence.


$$
\begin{array}{c}
  \\
  \\
  \\
  \\
1 \\
\mathbf{\tau_0}
\end{array} \qquad
\begin{array}{rcl}
 & & \\
 & & \\
 & & \\
 & 1 & \\
1\!\! & & \!\!1 \\
 &\mathbf{\tau_1} &
\end{array} \qquad
\begin{array}{rrcll}
 & & & & \\
 & & & & \\
 & & 1 & & \\
 & 2\!\! & & \!\!2 & \\
1\!\! & & 2 & & \!\!1 \\
 & & \mathbf{\tau_2} & &
\end{array} \qquad
\begin{array}{ccccccc}
 & & & & & & \\
 & & & 1 & & & \\
 & & 3\!\! & & \!\!3 & & \\
 & 3\!\! & & 6 & & \!\!3 & \\
1\!\! & & 3\!\! & & \!\!3 & & \!\!1 \\
 & & & \mathbf{\tau_3} & & &
\end{array} \qquad
\begin{array}{ccccccccc}
 & & & & 1 & & & & \\
 & & & 4\!\! & & \!\!4 & & & \\
 & & 6\!\! & & 12 & & \!\!6 & & \\
 & 4\!\! & & 12\!\! & & \!\!12 & & \!\!4 & \\
1\!\! & & 4\!\! & & 6 & & \!\!4 & & \!\!1 \\
 & & & & \mathbf{\tau_4} & & & &
\end{array}
$$

A {\it second multiplication principle} explains the pattern of $( \mathbf{\tau_n} )_n$:

\begin{lemma}\label{L:MP2}
Counting the rows from the top and the location on rows from left to right, the $j$-th number on the $i$-th row  in $\mathbf{\tau_n}$ equals
$\binom{n}{i-1} \cdot \binom{i-1}{j-1}$.
\end{lemma}

As with Lemma~\ref{L:MP}, we call this a {\it multiplication principle} not because the answer is given by the product of two binomial coefficients but because it can be viewed as a multiplication table of a sort, with rows indexed by the binomial coefficients and the other ``dimension" represented this time by the Pascal triangle. The proof is also by induction on $n$.

The triangle $\mathbf{\tau_n}$ gives the count of the number of random walks of length $n$ in $\R^3$ which start at the origin and end on the faces of the octahedron (the $l^1$-sphere of radius $n$ in $\R^3$). In terms of walks, the above second multiplication principle can be restated as follows: the number of random walks of length $n$ in the 3D-space which start at $(0,0,0)$ and end at point $(k,l,m)$, where $k, l, m\geq 0$ and $k+l+m=n$, is given by
$
\binom{n}{m} \cdot \binom{k+l}{l} = \binom{n}{k+l} \cdot \binom{k+l}{l}.
$
Perhaps surprisingly, it turns out that these numbers are exactly the ones that describe how {\it all} the paths distribute in the 3D-case, as long as we focus on how the distribution of paths happens in {\it entire planes} $z=h$, instead of at individual points. A few more definitions and notations are needed in order to justify this observation.

\begin{figure}[h]
\begin{center}
 \scalebox{.55}{\includegraphics{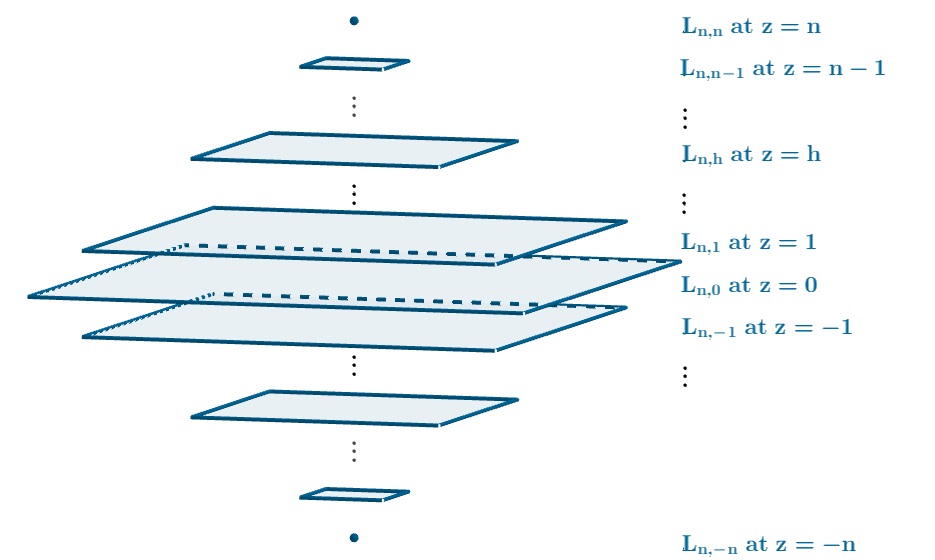}}
 \caption{Slicing the octahedron of possible end points with the ``layers" $L_{n,h}$.}
 \label{fig:slicing}
\end{center}
\end{figure}

For any integers $n\geq 1$ and $-n\leq h \leq n$, we call the {\it layer at height $h$ of walks of length $n$}, denoted $L_{n,h}$, the distribution of number of simple random walks of length $n$ that end in the plane $z=h$.
%
These layers provide a decomposition in horizontal ``slices" of the octahedron of paths as depicted in Figure~\ref{fig:slicing}.

\begin{figure}[h]
\begin{center}
\scalebox{.08}{\includegraphics{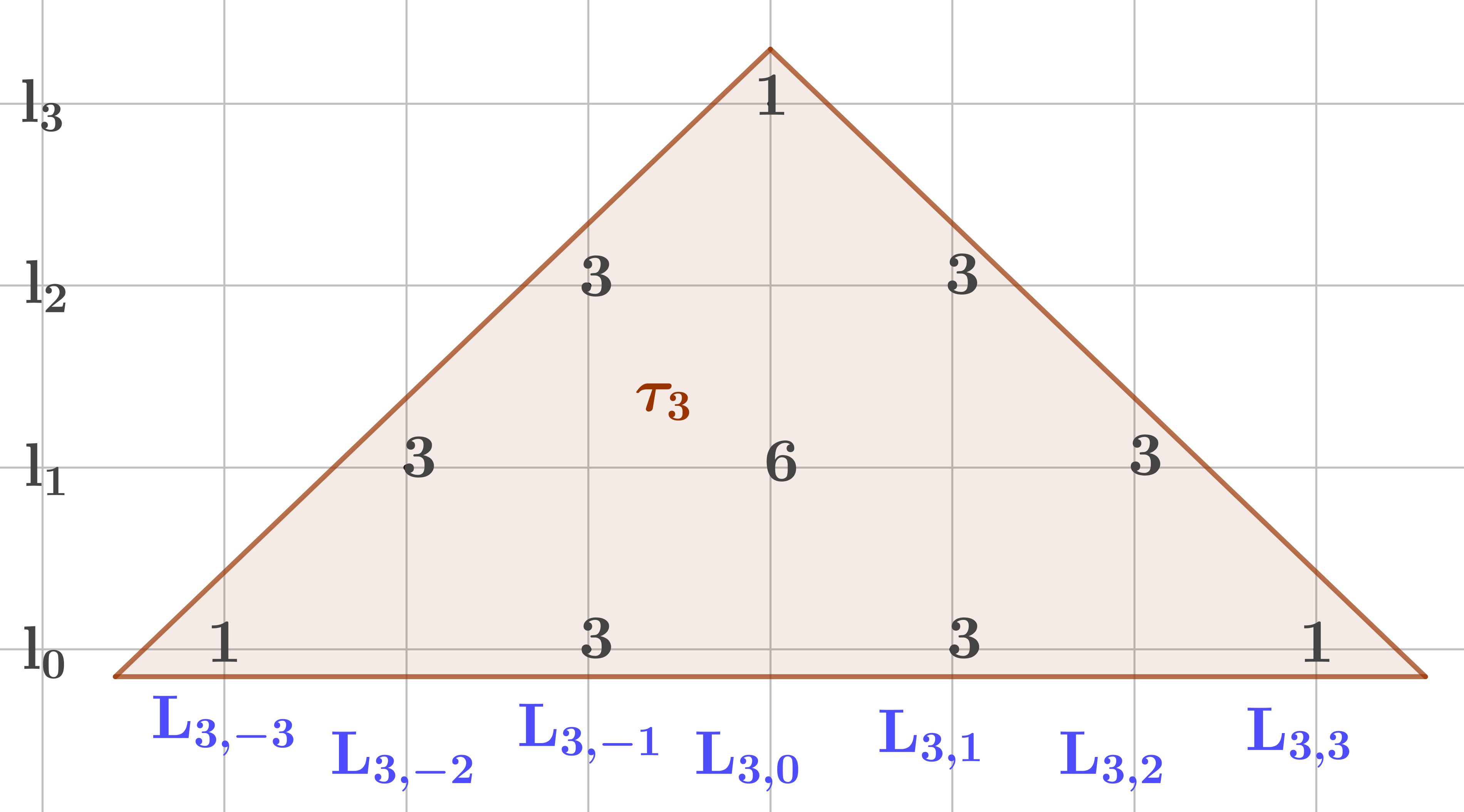}} \quad
\scalebox{.081}{\includegraphics{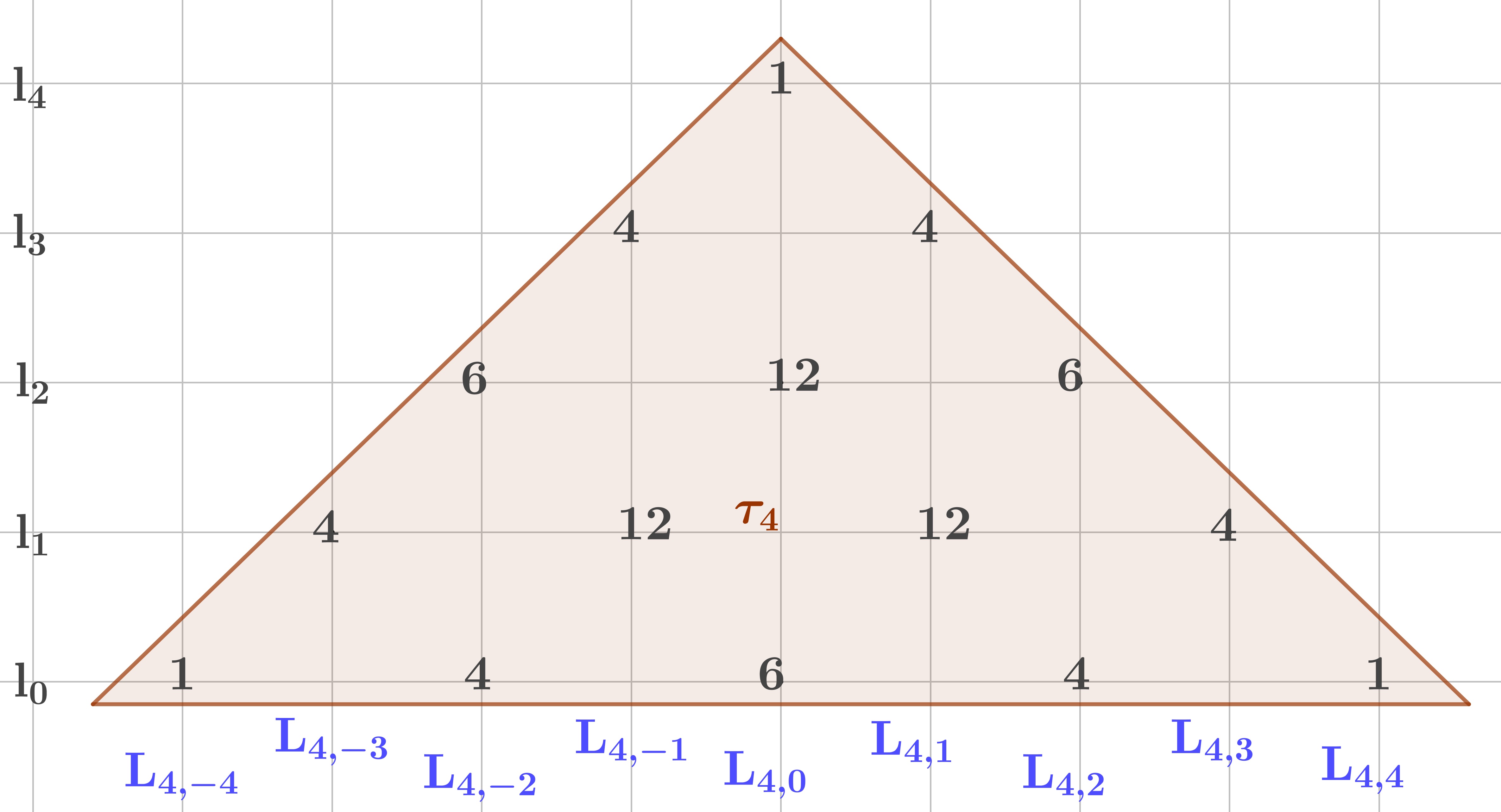}}
 \caption{Obtaining the $L_{n,h}$'s as linear combinations of the 2D-distributions with coefficients from $\mathbf{\tau_n}$, in the case when $n=3$ and $n=4$.}
 \label{fig:3d-layers}
\end{center}
\end{figure}

These layers can also be expressed in an exact form using the 2D-distribution of paths of length $n$, which we denote by $l_n$, provided by Lemma \ref{L:MP}. Imagine that the rows of the triangle $\mathbf{\tau_n}$ are formally multiplied with $l_n$, $l_{n-1}$, \dots, $l_1$, and $l_0$, respectively. Then, in this multiplied $\mathbf{\tau_n}$, the formal addition of the entries on each column gives the layers $L_{n,h}$, for $-n\leq h\leq n$. The addition is performed by keeping the 2D-distributions centered at the origin. Figure~\ref{fig:3d-layers} provides a visualization of this in the case when $n=3$ and $n=4$. For example, $L_{3,0}=6 l_1+l_3$ and $L_{4,0}=6 l_0+12 l_2 + l_4$.

Expanding the notation to any dimension $d\geq 1$, we denote by $l_n^{(d)}$ the distribution of {\it all} random walks of length $n$ in $\R^d$. We have:
$l_n^{(1)}$ is given by the binomial coefficients,
$l_n^{(2)}$ is a diamond, as described by Lemma~\ref{L:MP}, and
$l_n^{(3)}$ is an octahedron, as discussed above.
We also denote by $L_{n,h}^{(d+1)}$ the ``layer" of the distribution of random walks of length $n$ in $\R^{d+1}$ contained in the hyperplane $x_{d+1}=h$.
It is easy to see that $L_{n,n}^{(d+1)}=l_0^{(d)}$ and $L_{n,n-1}^{(d+1)}=n \cdot l_1^{(d)}$.
We formalize these observations in the following theorem, which shows how to obtain the geometric distribution of walks in dimension $d+1$ from the geometric distribution in dimension $d$, using the second multiplication principle.

\begin{theorem}\label{T:layers-allD}
With the above notations, for all $d\geq 1$, $n\geq 1$, and $0\leq h \leq n$,
\begin{equation} \label{eq:layers}
    L_{n,h}^{(d+1)} = \sum_{j=0}^{\lfloor \frac{n-h}{2} \rfloor}
        \binom{n}{n-h-2j} \binom{h+2j}{j} l_{n-h-2j}^{(d)}
        = \sum_{i=0 \text{ or } 1}^{n-h}
        \binom{n}{i} \binom{n-i}{\frac{n-h-i}{2}} l_{i}^{(d)}
\end{equation}
In the second sum, $i \equiv n-h \pmod{2}$.
\end{theorem}

\begin{proof}
Fix $d$. We will drop the superscript and keep in mind that any lower $l$ is the entire distribution of walks in $\R^d$ and any capital $L$ is just a layer in $\R^{d+1}$.
The proof is by induction on $n$. The case $n=1$ is clear: $L_{1,0}=l_1=\binom{1}{1} \binom{0}{0} l_1$ and $L_{1,1}=l_0=\binom{1}{0} \binom{1}{0} l_0$.

The key observation is that the random walks of length $n+1$ are obtained from the layers of the walks of length $n$, where each layer $L_{n,h}$ does only two ``moves:" either it {\it expands} horizontally in the hyperplane $x_{d+1}=h$, as if the walks were $d$-dimensional, or it {\it duplicates} above and below, in the hyperplanes $x_{d+1}=h+1$ and $x_{d+1}=h-1$, respectively.

\begin{figure}[h]
\begin{center}
\scalebox{.5}{\includegraphics{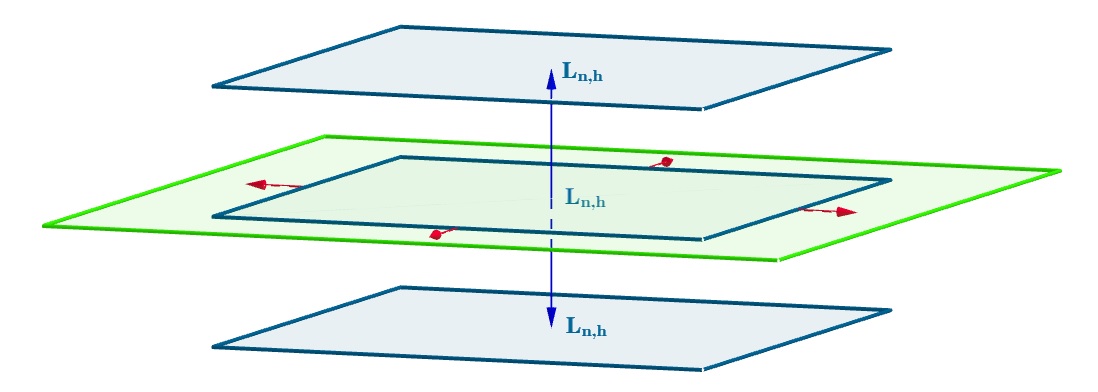}}
 \caption{Expansion and duplication of layers, as the length of paths increases.}
 \label{fig:3d-change-layers}
\end{center}
\end{figure}

\noindent If one interprets the $j$th column of $\mathbf{\tau_n}$ as keeping track of what is going on in the hyperplane $x_{d+1}=j-(n+1)$, then the observation that we have just made about the propagation of layers has the following interpretations. The horizontal expansion of a layer can be tracked by an upward addition ($l_k$ becomes $l_{k+1}$). The vertical duplications translate into left and right duplications (the horizontal hyperplanes are changed). But this is exactly how $\mathbf{\tau_{n+1}}$ is generated from $\mathbf{\tau_n}$. In other words, the sequence $(\mathbf{\tau_n})_n$ keeps track of the propagation of random walks!

\newpage
\begin{figure}[h]
\begin{center}
\scalebox{.15}{\includegraphics{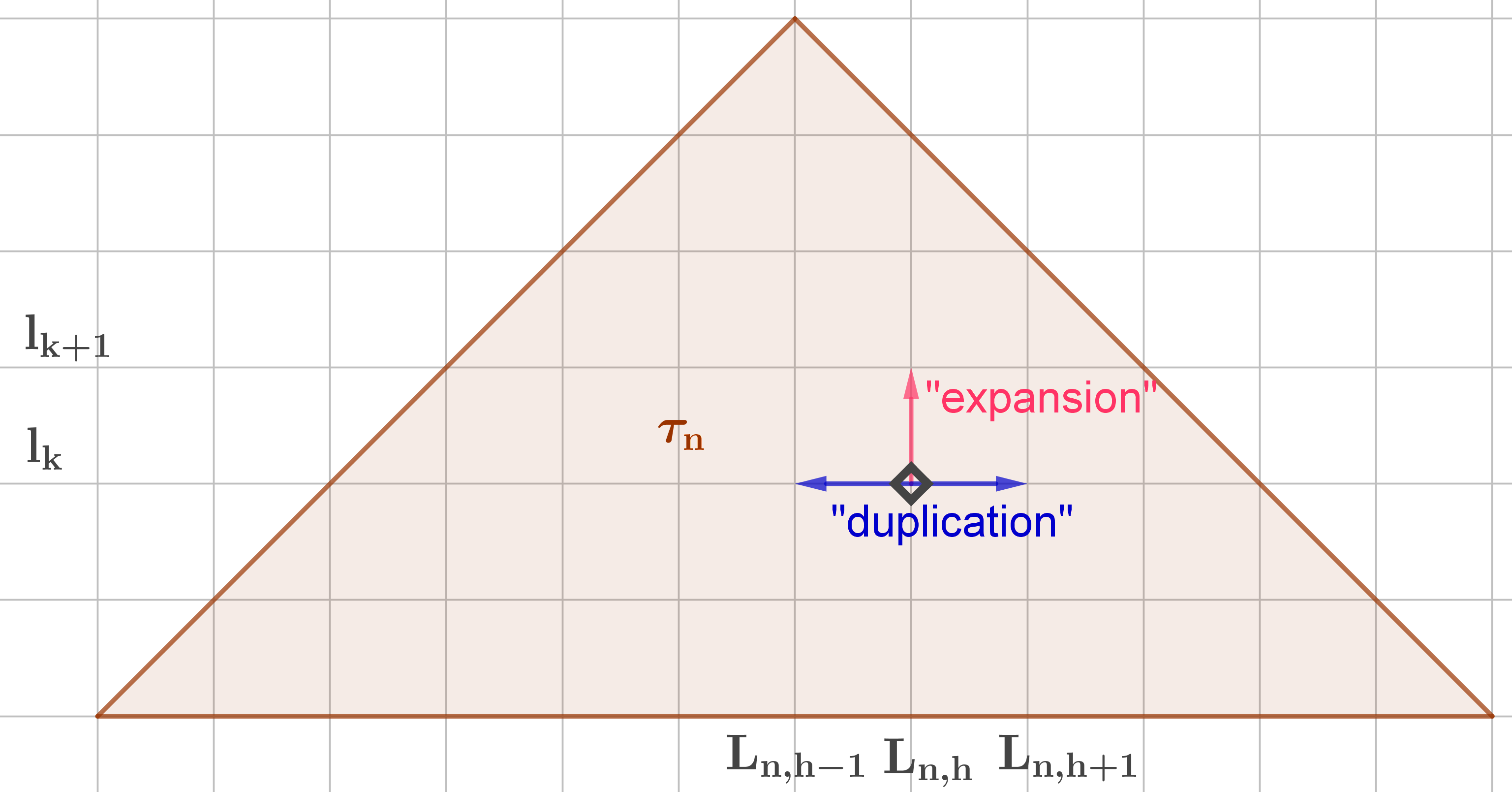}}
 \caption{Understanding the change in layers via $\mathbf{\tau_n}$.}
 \label{fig:2d-change-layers}
\end{center}
\end{figure}

\noindent Finally, the formula given in the statement of the theorem simply uses the numbers from $\mathbf{\tau_n}$ to compose $L_{n,h}$ out of the corresponding $d$-distributions $l_i$.
\end{proof}

For $m\geq 0$, we denote by $A_{2m}^{(d)}$ the number of walks of length $2m$ in $\R^d$ that return at the origin at time $t=2m$. In the literature these are known as {\it closed walks} or {\it excursions}. We will use the first term in our paper. By convention $A_{0}^{(d)}=1$.

\begin{cor} \label{C:As-all-d}
For any $d\geq 1$,

\begin{equation} \label{eq:As-all-d} 
A_{2m}^{(d+1)} = \binom{2m}{m} \cdot \sum_{k=0}^{m} \binom{m}{k}^2 \,
        \frac{A_{2k}^{(d)}}{\binom{2k}{k}}.
\end{equation}

\end{cor}

\begin{proof}
Use Theorem~\ref{T:layers-allD} for $n=2m$ to evaluate the layer corresponding to $h=0$ at the origin:
{\small
$$
\begin{aligned}
A_{2m}^{(d+1)} &= \sum_{k=0}^{m} \binom{2m}{2k} \binom{2m-2k}{m-k} A_{2k}^{(d)}
   = \sum_{k=0}^{m}
    \frac{(2m)! \cdot (2m-2k)!}{(2k)!\cdot(2m-2k)!\cdot(m-k)!^2} \, A_{2k}^{(d)} \\
   & = \binom{2m}{m} \cdot \sum_{k=0}^{m}
    \frac{m!^2}{k!^2 \cdot(m-k)!^2} \, \frac{A_{2k}^{(d)}}{\binom{2k}{k}}
     = \binom{2m}{m} \cdot \sum_{k=0}^{m} \binom{m}{k}^2 \,
        \frac{A_{2k}^{(d)}}{\binom{2k}{k}}.
\end{aligned}
$$
}
\end{proof}

\section{Asymptotics of the number of closed walks}
In this section we study the asymptotic behaviour of the number of closed walks (with starting point assumed for simplicity to be the origin).
For any $d\geq 1$, recall that, before Corollary~\ref{C:As-all-d}, we denoted by $A_{2n}^{(d)}$ the number of closed walks of length $2n$ in $\R^d$. Similarly, we denote by $B_{2n}^{(d)}$ the number of walks of length $2n$ in $\R^d$ that return {\it for the first time} at the origin at time $t=2n$ (that is, which avoided the origin till the very end).

In the one-dimensional case, we have, of course,
$A_{2n}^{(1)}= \binom{2n}{n}$
(\href{https://oeis.org/search?q=A000984}{A000984})
and $B_{2n}^{(1)}$ is the sequence
\href{https://oeis.org/search?q=A284016}{A284016}
in \cite{Sloane}. The just cited reference from the OEIS gives the formula $B_{2n}^{(1)}=A_{2n}^{(1)}/(2n-1)$, which we will justify further below (see Remark \ref{r:b-formula1D}).
In the two-dimensional case, we have, based on Lemma~\ref{L:MP}, $A_{2n}^{(2)}= \binom{2n}{n}^2$
(\href{https://oeis.org/search?q=A002894}{A002894})
and $B_{2n}^{(2)}$ is the sequence
\href{https://oeis.org/search?q=A054474}{A054474}.
There seems to be no formula involving combinatorics for $B_{2n}^{(2)}$'s.
In the three-dimensional case, Corollary~\ref{C:As-all-d} gives us formulas for the $A_{2n}^{(3)}$'s (\href{https://oeis.org/search?q=A002896}{A002896})
and $B_{2n}^{(3)}$ is the sequence
\href{https://oeis.org/search?q=A049037}{A049037}.
In the four-dimensional case, Corollary~\ref{C:As-all-d} gives us formulas for the $A_{2n}^{(4)}$'s (\href{https://oeis.org/search?q=A039699}{A039699}).
At the time our preprint \cite{DS} was posted, the sequence of $B_{2n}^{(4)}$'s was not found in \cite{Sloane}. Now it is \href{https://oeis.org/A359801}{A359801}.

The first eight terms of these sequences are provided in the table below.

\begin{center}
\begin{tabular}{|c||c|c|c|c|c|c|c|c|}
\hline
$n$ & 1 & 2 & 3 & 4 & 5 & 6 & 7 & 8 \\
\hline \hline
$A_{2n}^{(1)}$ & 2 & 6 & 20 & 70 & 252 & 924 & 3,432 & 12,870 \\
\hline
$B_{2n}^{(1)}$ & 2 & 2 & 4 & 10 & 28 & 84 & 264 & 858 \\
\hline \hline
$A_{2n}^{(2)}$ & 4 & 36 & 400 & 4,900 & 63,504 & 853,776 & 11,778,624 & 165,636,900 \\
\hline
$B_{2n}^{(2)}$ & 4 & 20 & 176 & 1,876 & 22,064 & 275,568 & 3,584,064 & 47,995,476 \\
\hline \hline
$A_{2n}^{(3)}$ & 6 & 90 & 1,860 & 44,730 & 1,172,556 & 32,496,156 & 936,369,720 & 27,770,358,330 \\
\hline
$B_{2n}^{(3)}$ & 6 & 54 & 996 & 22,734 & 577,692 & 15,680,628 & 445,162,392 & 13,055,851,998 \\
\hline \hline
$A_{2n}^{(4)}$ & 8& 168& 5,120& 190,120& 7,939,008& 357,713,664& 16,993,726,464& 839,358,285,480 \\
\hline
$B_{2n}^{(4)}$ & 8& 104& 2,944& 108,136& 4,525,888& 204,981,888& 9,792,786,432& 486,323,201,640 \\
\hline \hline
\end{tabular}
\end{center}

\noindent It is interesting to note that when $d=3$ and $d=4$ the $B_{2n}^{(d)}$'s are fairly comparable with the $A_{2n}^{(d)}$'s, as opposed with the 1D- and 2D-case, where the first returns are much smaller than the returns. We will explain this behaviour and more in this section.

\subsection{Recurrence relation for the first time returns}\label{subsec:recurrence}
We show first how to generate recursively the sequences $(B_{2n}^{(d)})_n$ from the corresponding sequences $(A_{2n}^{(d)})_n$.
Because in this recurrence relation the dimension does not play a role, we simplify temporarily the notation and write $\alpha_n=A_{2n}^{(d)}$ and $\beta_n=B_{2n}^{(d)}$. It is clear that we have $\beta_1=\alpha_1$, as the walks that return in two steps at the origin return there for the first time as well.
One then observes that $\beta_2=\alpha_2 - \beta_1 \alpha_1$, because from the walks that return at the origin after 4 steps (namely, $\alpha_2$) we have to subtract those generated after 2 steps by the $\beta_1$ walks that have already returned (and there are $\beta_1 \alpha_1$ such walks in total).
This argument generalizes to give, for all $n\geq 2$,
\begin{equation}\label{eq:recurrenceBs}
  \beta_n = \alpha_n -\sum_{k=1}^{n-1} \beta_k \alpha_{n-k},
  \text{ or }
  B_{2n}^{(d)} = A_{2n}^{(d)} -\sum_{k=1}^{n-1} B_{2k}^{(d)} A_{2n-2k}^{(d)}.
\end{equation}

\noindent
The identity~(\ref{eq:recurrenceBs}) can be written as the following identity of series:
\begin{equation}\label{eq:ab-series}
  \left( 1 - \beta_1 x - \beta_2 x^2 - \cdots - \beta_n x^n - \cdots \right) \cdot
    \left( 1 + \alpha_1 x + \alpha_2 x^2 + \cdots +\alpha_n x^n + \cdots \right) = 1.
\end{equation}

The identity~(\ref{eq:ab-series}) has two main consequences. First, once we establish the asymptotics of the $A_{2n}^{(d)}$'s, it will be used to obtain the asymptotics of the $B_{2n}^{(d)}$'s. See Subsection~\ref{subsec:Bs}. Second, it provides a proof of Gy\"{o}rgy P\'{o}lya's celebrated theorem about the recurrence and transience of random walks (\cite{Polya}). We recall that an infinite simple random walk is called {\it recurrent} if it is certain to return at its starting point; if not, the random walk is called {\it transient}. P\'{o}lya's Theorem asserts that an infinite simple random walk on a $d$-dimensional lattice is recurrent for $d=1$ and $d=2$, and it is transient for $d\geq 3$.
Theorem~\ref{T:TheoremAs} below implies that, when $d \ge 3$, the series of $A_d(w)$ converges absolutely for $|w|\le \frac{1}{(2d)^2}$ and so does the series of $B_d(w)$, since $0< B_{2n}^{(d)} \le A_{2n}^{(d)}$. Hence the relation $(1-B_d(w)) A_d(w)=1$ extends to the circle of convergence by continuity, so $1-B_d(1/(2d)^2)=1/m_d$, with $0<1/m_d<1$. It follows that $0 < p_d = B_d(1/(2d)^2) < 1$.

The proof of the recurrence of the random walk when $d=1$ and $d=2$ follows directly.
When $d=1$, recalling that $\alpha_n=A_{2n}^{(1)}= \binom{2n}{n}$, the series with the $\alpha_n$'s is identified using Newton's generalized binomial formula for $r=-1/2$:
%
%
$$
1 + \sum_{n=1}^{\infty} \alpha_n x^n = (1-4x)^{-1/2}.
$$
This implies that
\begin{equation}\label{eq:b-series1D}
  1 - \beta_1 x - \beta_2 x^2 - \dots -\beta_n x^n - \dots = (1-4x)^{1/2},
\end{equation}
with interval of convergence $[-1/4,1/4]$. The probability of return to the origin, namely $\sum_{n=1}^\infty \beta_n/2^{2n}$, equals 1 by evaluating (\ref{eq:b-series1D}) at $x=1/4$.

\vspace{.3cm}
\begin{remark}\label{r:b-formula1D}
From Equation~(\ref{eq:b-series1D}), it is now clear that $\beta_n=B_{2n}^{(1)}= A_{2n}^{(1)}/(2n-1)$.
\end{remark}

\vspace{.3cm}
The two-dimensional case is less elementary. 
Recalling that $\alpha_n=A_{2n}^{(2)}= \binom{2n}{n}^2$, these coefficients are recognized in the expansion of the complete elliptic integral of the first kind:
%
%
$$
\frac{K(k)}{\pi/2} = \tensor[_2]{F}{_1} \left( \tfrac{1}{2}, \tfrac{1}{2}; 1; k^2 \right)
    = \sum_{i=0}^{\infty} \binom{2i}{i}^2 \frac{1}{16^i} \cdot k^{2i},
$$
which holds true for $|k|<1$.
This implies that
\begin{equation}\label{eq:b-series2D}
1 - \beta_1 x - \beta_2 x^2 - \dots - \beta_n x^n - \dots
  = \frac{1}{\tensor[_2]{F}{_1}\left( \frac12, \frac12; 1; 16 x \right)},
\end{equation}
with interval of convergence $(-1/16,1/16)$.
The limit when $x\rightarrow 1/16^{-}$ of $\tensor[_2]{F}{_1}$ is infinite and
it follows again that the probability of return to the origin equals 1, since the $\beta_n$'s are positive so we can move the limit inside.

\subsection{Asymptotics of the numbers $A_{2n}^{(d)}$'s of closed walks }\label{subsec:As}

For $d\geq 1$, we denote
$\displaystyle{ x_n^{(d)}=\frac{A_{2n}^{(d)}}{{2n \choose n}} }$,
with $x_0^{(d)}=1$.\footnote{In \cite{Sloane},
$(x_n^{(3)})_n$ is
\href{https://oeis.org/search?q=A002893}{A002893}
and $(x_n^{(4)})_n$ is
\href{https://oeis.org/search?q=A002895}{A002895}.
}
The generating formula for the $A_{2n}^{(d)}$'s contained in Equation (\ref{eq:As-all-d}) of Corollary~\ref{C:As-all-d}, gives the following fundamental recurrence relation:

\begin{equation} \label{eq:Xs} 
x_n^{(d+1)} = \sum_{k=0}^n {n \choose k}^2 x_k^{(d)} =
   1+\sum_{k=1}^n{n \choose k}^2 x_k^{(d)}, \text{ for } n \ge 1, d \ge 1.
\end{equation}

\begin{theorem}\label{T:TheoremA}
The following hold:
\begin{enumerate}
\item[(i)] $$
x_{n}^{(d)}=a_d\frac{d^{2n}}{(\pi n)^{(d-1)/2}} \left(1+ \frac{r_1(d)}{n}+\frac{r_2(d)}{n^2}+\dots +\frac{r_m(d)}{n^m}+
O_{d,m}\left( \frac{1}{n^{m+1}} \right) \right).
$$

where $\displaystyle{ a_d=\frac{d^{d/2}}{2^{d-1}} }, r_1(d)=\frac{1-d}{8}, r_2(d)=\frac{(d^2-1)(2d-3)}{384}$, while $r_3(d), r_4(d)$ are polynomials in $d$ of degree $4$ and $6$ (see below), and we conjecture that all $r_m(d)$ are polynomials in $d$ of degree $[\frac{3m}{2}]$ and divisible by $d-1$.

\item[(ii)]
$$
A_{2n}^{(d)}=a_d\frac{(2d)^{2n}}{(\pi n)^{d/2}} \left(1+ \frac{a_1(d)}{n}+\frac{a_2(d)}{n^2}+\dots +\frac{a_m(d)}{n^m}+
O_{d,m}\left( \frac{1}{n^{m+1}} \right) \right).
$$

where $\displaystyle{ a_1(d)= -\frac {d}{8}, a_2(d)=\frac{2d^3-3d^2+4d}{384}}$, while $a_3(d), a_4(d)$ are polynomials in $d$ of degree $4$ and $6$, and we conjecture that all $a_m(d)$ are polynomials in $d$ of the same degree as $r_m(d)$.

\end{enumerate}
\end{theorem}

\vspace{.3cm}
\begin{remark} The polynomials mentioned in the theorem are:

\begin{equation}
 \begin{aligned}
  r_3(d) = & \tfrac{1}{3072} \left( d-1 \right)  \left( 6\,{d}^{3}-19\,{d}^{2}+14\,d+15 \right), \\
  r_4(d) = & {\tfrac {1}{1474560}}\, \left( d-1 \right)  \left( 20\,{d}^{5}+2504\,{d
    }^{4}-10241\,{d}^{3}+9679\,{d}^{2}+309\,d+945 \right); \\
  a_3(d) = & {\tfrac {1}{1024}}\,{d}^{2} \left( 2\,{d}^{2}-9\,d+12 \right), \\
  a_4(d) = & {\tfrac {1}{1474560}}\,d \left( 20\,{d}^{5}+2484\,{d}^{4}-13105\,{d}^{3
    }+21480\,{d}^{2}-11440\,d-384 \right).
 \end{aligned}
\end{equation}

\end{remark}

\vspace{.3cm}
\begin{remark}
Using the well known asymptotic expansion
${2n \choose n} \sim \frac{2^{2n}}{\sqrt {\pi n}} \left( 1+\frac{g_1}{n}+\cdots+\frac{g_k}{n^k}+\cdots \right)$, where $g_1=-\frac{1}{8}, g_2 =\frac{1}{128}, g_3=\frac{5}{1024}, g_4=-\frac{21}{32768}$, and so on, then clearly $a_m(d)=r_m(d)+g_1r_{m-1}(d)+\dots +g_m$, so if $r_k(d)$, $k=1, \dots, m$, are polynomials, so are $a_k(d)$, $k=1, \dots, m$, and viceversa. Moreover, if the degree conjecture is correct (as supported by our Legendre polynomials framework, as we will show in the proof of Theorem~\ref{T:TheoremA}), then $a_m(d)$ and $r_m(d)$ have the same leading coefficient.
\end{remark}

\vspace{.3cm}
\begin{remark}
Because $x_n^{(1)}=1$ it follows that $r_m(1)=0$, $m \ge 1$, and similarly $r_m(2)=g_m$, $m \ge 1$. The values of $a_1(d)$, $a_2(d)$ appear in \cite{Domb} and a few values of $a_m(3)$ are given in \cite{Joyce73}. Our computations are completely independent of these, as they use the new Legendre polynomials framework with which we prove Theorem~\ref{T:TheoremA}, and are consistent with the results existing in the literature.
\end{remark}

The proof by induction on $d \ge 2$ will be elementary and self-contained, using only basic calculus and classical results about Legendre polynomials, which we recall first. We note that a similar formula for the $a_d$'s appears indirectly in the original paper \cite[eqn.(12)]{Polya} and these asymptotics are developed in \cite{Domb} using P\'{o}lya's original method, which is completely different than ours.

We recall (\cite[15.1-9]{WW}, \cite[4.8]{Szego}) that the Legendre polynomials $P_n$ of degree $n \ge 0$, which form an orthogonal basis in the space of real polynomials on $[-1,1]$ with the usual $L^2$-norm, are given by Rodrigues' formula
$\displaystyle{
P_n(x)=\frac{1}{2^n n!}\frac{d^n}{dx^n}(x^2-1)^n.
}$
The first ones are $P_0(x)=1, P_1(x)=x, P_2(x)=\frac{3x^2-1}{2}$, and so on. By Leibniz rule, the above leads to the formula:

\begin{equation} \label{eq:Legendre_recurrence}
P_n(x)=\left( \frac{x-1}{2} \right)^n \sum_{k=0}^n {n \choose k}^2
  \left( \frac{x+1}{x-1} \right)^k.
\end{equation}

\noindent
There is also an integral formula
$$
P_n(x)=\frac{1}{2\pi}\int_{-\pi}^{\pi}(x+\sqrt {x^2-1}\cos \theta)^n \, d\theta,
$$
valid for all complex $x$, as it is clear that the choice of the square root is irrelevant as all odd powers in the binomial expansion vanish.

For any (complex number) $x \notin [-1,1]$ there is a unique complex number $z$ satisfying $|z|>1$ and $x=\frac{1}{2}(z+\frac{1}{z})$. This is because the quadratic in $z$ solving for $x$ is symmetric under $z \to 1/z$, and if $z=e^{i\theta}$ then $x=\cos \theta \in [-1,1]$. Note that then $\Re (1-z^{-2})>0$, so one has natural definition of $(1-z^{-2})^a$ for any $a \in \mathbb R$ which is positive for $z>1$.
In particular, for $x>1, z=x+\sqrt {x^2-1}>1$ and the integral representation above gives

\begin{equation} \label{eq:1star} 
|P_n(x)| \le |z|^n.
\end{equation}

The fundamental asymptotic theorem (Laplace-Heine) for the Legendre polynomials (of argument not in $[-1,1]$) states the following (\cite[8.21]{Szego}). Given a fixed Jordan curve $\Gamma$ enclosing the interval $[-1,1]$ one has uniformly for $x$ in the exterior of $\Gamma$:

$$
P_n(x) = q_n z^n \sum_{j=0}^{m-1} q_j \frac{1 \cdot 3 \cdots (2j-1)}{(2n-1)(2n-3) \cdots (2n-2j+1)} z^{-2j} (1-z^{-2})^{-j-1/2} + O_{\Gamma,m}\left( \frac{|z|^n}{n^{1/2+m}} \right),
$$
for all $m \ge 1$, where $q_0=1$ and $q_k={2k \choose k}/2^{2k}$, for $k\geq 1$.
In the above $x=\frac{1}{2}(z+\frac{1}{z}), |z|>1$, so $(1-z^{-2})^{-1/2-k}$ are naturally well defined as noted.

Let now $x, y>1$ satisfy $x=\frac{y+1}{y-1}$, so that $z=\frac{(\sqrt y+1)^2}{y-1}$ and $(1-z^{-2})^{-1/2}=\frac{\sqrt y+1}{2y^{1/4}}$. Formula~(\ref{eq:Legendre_recurrence}) becomes

\begin{equation}
\label{eq:Legendre_recurrence2} 
(y-1)^n P_n(x)=\sum_{k=0}^n {n \choose k}^2 y^k,
\end{equation}

\noindent
while the asymptotic above becomes

\begin{equation} \label{eq:Legendre_asymptotics}
\sum_{k=0}^n {n \choose k}^2 y^k = q_n\frac{(\sqrt y+1)^{2n+1}}{2y^{1/4}}
\sum_{j=0}^{m-1} q_{jn} \frac{(\sqrt y-1)^{2j}}{4^j y^{j/2}} +
O_{\Gamma,m}\left( \frac{(\sqrt y+1)^{2n}}{n^{1/2+m}} \right),
\end{equation}
where $q_{0n}=q_0=1$ and $q_{jn}=q_j \frac{1 \cdot 3 \cdots (2j-1)}{(2n-1)(2n-3) \cdots (2n-2j+1)}$. The growth of the coefficients appearing in equation~(\ref{eq:Legendre_asymptotics}) is worth noticing: $q_n= {2n \choose n}/{2^{2n}} \sim {1}/{\sqrt {\pi n}}$ and $q_{jn} \sim c_j n^{-j}$.
We also remark that in the left-hand side of the sum of this equation we can ignore the first term, for $k=0$ which is $1$, because it is negligible with respect the right-hand side error.

With the same notation, the inequality in (\ref{eq:1star}) translates to:

\begin{equation} \label{eq:4stars} 
\sum_{k=0}^n {n \choose k}^2 y^k \le (\sqrt y+1)^{2n}.
\end{equation}

We start the proof of Theorem~\ref{T:TheoremA} by first showing that in the induction step the asymptotic carryover goes through relation (\ref{eq:Xs}). Using a simple version of the elimination method in asymptotic analysis, we reduce our induction step to developing asymptotics for integrals where it is straightforward to apply Laplace method as described in \cite[Appendix~B.6]{FS}. We perform the standard steps of localization to the smaller range $\frac{\log^2 n}{n}$ (``Neglect the tails"), approximate the integrand by forcing out the exponential main term $(d+1)^{2n+1}$  and using Taylor series on the resulting function (``Centrally approximate"), and finally introduce back the tails (``Complete the tails").
First, in Lemma \ref{L:carryover} we deal with the asymptotic error carryover. In its statement and proof $C_{n,m,d}$, $C_{m,d}$,  and so on are constants with the implied dependencies, but can change between uses.

\begin{lemma}\label{L:carryover}
Assume that for $d \ge 2, n \ge 1, m\ge 0$ we prove that
$$
x_{n}^{(d)}=a_d\frac{d^{2n}}{(\pi n)^{(d-1)/2}} \left(1+ \frac{r_1(d)}{n}+\frac{r_2(d)}{n^2}+\dots +\frac{r_m(d)}{n^m}+
 \frac{C_{n,m,d}}{n^{m+1}} \right),
$$
with $|C_{n,m,d}| \le C_{m,d}$. Then it follows that
$$
x_n^{(d+1)}=a_d\sum_{k=1}^n{n \choose k}^2\frac{d^{2k}}{(\pi k)^{(d-1)/2}} \left(1+ \frac{r_1(d)}{k}+\frac{r_2(d)}{k^2}+\dots +\frac{r_m(d)}{k^m}\right) +C_{n,m,d+1}\frac{(d+1)^{2n}}{n^{d/2+m+1}},
$$
with $|C_{n,m,d+1}| \le C_{m,d+1}$.

\end{lemma}

\begin{proof}
We substitute the above expression for $x_{k}^{(d)}, k \ge 1$, in formula (\ref{eq:Xs}) and note that the first term $1$ is clearly subsumed in the final error. For $0 \le k \le n/4$, using that ${n \choose k}^2 \le 4^n$, we have ${n \choose k}^2{d^{2k}} \le (4\sqrt d)^n$ so
$$
\left| \sum_{k \le n/4}{n \choose k}^2\frac{d^{2k}}{(\pi k)^{(d-1)/2}} \frac{C_{k,m,d}}{k^{m+1}}
    \right| \le C_{m,d} \cdot n \cdot (4\sqrt d)^n =
    O_{m,d}\left( \frac{(d+1)^{2n}}{n^{d/2+m+1}} \right),
$$
since $\frac{4\sqrt d}{(d+1)^2}<1$.

For $k >n/4$ we have $|C_{k,m,d}|k^{(1-d)/2-m-1} \le  4^{(d-1)/2+m+1}C_{m,d}n^{(1-d)/2-m-1}$, so the absolute value of this other part of the sum that gives the error terms is at most (with a different $C_{m,d}$)
$$
C_{m,d} \frac{1}{ n^{(d-1)/2+m+1}} \cdot \sum_{k > n/4}{n \choose k}^2 d^{2k}
    \le C_{m,d+1}\frac{(d+1)^{2n}}{n^{d/2+m+1}}.
$$
The last inequality comes from (\ref{eq:Legendre_asymptotics}) with $y=d^2$, remembering that in formula (\ref{eq:Legendre_asymptotics}) $q_n=O(1/\sqrt n)$ and  the inner sum is $O_{d,m}(1)$.
\end{proof}

The following two lemmas are used in justifying the steps in the Laplace method indicated above.

\begin{lemma}\label{L:integral-estimate}
If $a,b>0, k \ge 1$ then

$$\int_b^{\infty}e^{-sk}s^{a-1}ds \le \frac{(b+1)^a(a+1)^{a+1}}{\min (b,1)}e^{-bk}.$$
\end{lemma}

\begin{proof}
Use induction on the integer part of $a$ and integration by parts.
\end{proof}

\begin{lemma}\label{L:exponential-estimate}
If $0 \le q \le 1/2$ then $(e^{-q} d+1)^{2n} \le (d+1)^{2n}e^{-qnd/(d+1)}$, for $d \ge 2$ and $n \ge 1$.
\end{lemma}

\begin{proof}
We use, for $0\leq x<1$,  that $\log(1-x) \le -x$ and $1-e^{-x}\ge x/2$. Then
$(e^{-q} d+1)^{2n}=(d+1)^{2n}(1-\frac{(1-e^{-q})d}{d+1})^{2n}$,
hence taking logarithms
$2n\log (1-\frac{(1-e^{-q})d}{d+1}) \le -2n \frac{(1-e^{-q})d}{d+1} \le -\frac{qnd}{d+1}$.
\newline
\end{proof}

We use now a simple version of the elimination method from \cite[3.9]{Dingle}. From

$$
k^{(1-d)/2-j}=\frac{1}{\Gamma( (d-1)/2+j )}
\int_0^{\infty}e^{-sk}\cdot s^{(d-1)/2+j-1} \, ds,
\text{ for } k \ge 1, d \ge 2, j \ge 0,
$$
we can write (with $r_0(d)=1$ and noting that since $d \ge 2$ we have integrability near $0$)

\begingroup\makeatletter\def\f@size{11}\check@mathfonts
$$
\sum_{k=1}^n {n \choose k}^2 \frac{d^{2k} \cdot r_j(d)}{(\pi k)^{(d-1)/2} \cdot k^j}
    =\frac{r_j(d)}{\pi^{(d-1)/2}\cdot \Gamma( (d-1)/2+j )} \int_0^{\infty}
     \sum_{k=1}^n {n \choose k}^2 (e^{-s}d^2)^k \cdot s^{(d-1)/2+j-1} \,ds.
$$
\endgroup

\noindent
So, under the assumptions of Lemma \ref{L:carryover}, we obtain

\begingroup\makeatletter\def\f@size{10}\check@mathfonts
\begin{equation} \label{eq:xns}
x_n^{(d+1)} = a_d \sum_{j=0}^m
  \frac{r_j(d)}{\pi^{(d-1)/2} \cdot \Gamma( (d-1)/2+j )}
    \int_0^{\infty} \sum_{k=1}^n{n \choose k}^2 (e^{-s}d^2)^k \cdot s^{(d-1)/2+j-1} \,ds
       + O_{m,d}\left( \frac{(d+1)^{2n}}{n^{d/2+m+1}} \right).
\end{equation}
\endgroup

\noindent
Applying Lemma \ref{L:integral-estimate} with $b_n=\frac{\log^2 n}{n}$, $n \ge 2$, $b_n >0$, $k \ge 1$, inequality (\ref{eq:4stars}), and Lemma~\ref{L:exponential-estimate}, we see that we can neglect the tails in (\ref{eq:xns}) up to allowable error, so we get

\begingroup\makeatletter\def\f@size{10}\check@mathfonts
\begin{equation} \label{eq:xns-expanded}
x_n^{(d+1)}=a_d\sum_{j=0}^m\frac{r_j(d)}{\pi^{(d-1)/2}\Gamma((d-1)/2+j)}\int_0^{b_n}
     \sum_{k=1}^n{n \choose k}^2(e^{-s}d^2)^k \cdot s^{(d-1)/2+j-1} \,ds +
     O_{m,d}\left( \frac{(d+1)^{2n}}{n^{d/2+m+1}} \right).
\end{equation}
\endgroup

\noindent
Since $e^{-b_n}d^2 \sim d^2 \geq 4$, for $d \ge 2$, we can substitute (\ref{eq:Legendre_asymptotics}) with $y=e^{-s}d^2$, $0 \le s \le b_n$, in the above.
Putting $t=ns$ and using Lemma \ref{L:exponential-estimate} for $j=0, \dots ,m$ we estimate the error coming from (\ref{eq:Legendre_asymptotics}) as follows:

\begingroup\makeatletter\def\f@size{11}\check@mathfonts
$$
\int_0^{b_n}(e^{-s/2}d+1)^{2n}s^{(d-1)/2+j-1} \,ds \le \frac{(d+1)^{2n}}{n^{(d-1)/2+j}}\int_0^{\log^2 n}e^{-td/(2d+2)} t^{(d-1)/2+j-1} \, dt <<_{d,m}\frac{(d+1)^{2n}}{n^{(d-1)/2+j}}.
$$
\endgroup

\noindent
So the error term coming from (\ref{eq:Legendre_asymptotics}) is $<<_{d,m}\frac{(d+1)^{2n}}{n^{(d-1)/2+j+1/2+m}}<<_{d,m}\frac{(d+1)^{2n}}{n^{d/2+m+1}}$ as required. Note that the above considerations and the fact that $q_{kn} \sim c_k n^{-k}$ show that, if we fix $j=0, \dots ,m$ when we apply the asymptotic (\ref{eq:Legendre_asymptotics}), we actually need to go only up to $k=m-j$ terms as then the next term is $<<_{d,m}\frac{(d+1)^{2n}}{n^{d/2+m+1}}$ so is included in the required error.

Hence up to allowable error $O_{m,d}\left( \frac{(d+1)^{2n}}{n^{d/2+m+1}} \right)$, we have

\begin{equation} \label{eq:xns-expanded2}
x_n^{(d+1)} = \frac{a_d \cdot q_n}{2 \sqrt{d} \cdot \pi^{(d-1)/2}}
    \sum_{j=0}^m \sum_{k=0}^{m-j} \frac{r_j(d) \cdot q_{kn}}{(4d)^k \cdot \Gamma((d-1)/2+j)}
        I_{djk},
\end{equation}
where
\begin{equation} \label{eq:Idjk}
I_{djk} = \int_0^{b_n}
     (e^{-s/2}d+1)^{2n+1} \cdot (e^{-s/2}d-1)^{2k} \cdot
        e^{s/4+sk/2} \cdot s^{(d-1)/2+j-1} \,ds.
\end{equation}

\noindent
Putting $t=ns$ and forcing $(d+1)^{2n+1}$ out of the integral gives us

\begin{equation} \label{eq:xns-expanded3}
x_n^{(d+1)}= a_{n,d}\sum_{j=0}^m \sum_{k=0}^{m-j} \tilde{q}_{jk} \tilde{I}_{djk},
\end{equation}
where
\begin{equation} \label{eq:Idjk-tilde}
\tilde{I}_{djk} = \int_0^{\log^2 n} \left( 1-(1-e^{-t/(2n)}) \tfrac{d}{d+1} \right)^{2n+1}
    \cdot (e^{-t/(2n)}d-1)^{2k}
        \cdot e^{t/(4n)+tk/(2n)}
            \cdot t^{(d-1)/2+j-1} \,dt,
\end{equation}
and
\begin{equation} \label{eq:xns-expanded3-coeffs}
a_{n,d}=\frac{a_d \cdot q_n \cdot (d+1)^{2n+1}}{2\sqrt{d} \cdot (\pi n)^{(d-1)/2}}
\text{ and }
\tilde q_{jk}=\frac{r_j(d) \cdot q_{kn}}{(4d)^k \cdot \Gamma((d-1)/2+j)\cdot n^j}
    \sim c_{d,j,k} \cdot n^{-j-k}.
\end{equation}
(In the equation above, $\tilde q_{jk}$ has an asymptotic in inverse powers of $n$ starting at $n^{-j-k}$ since $q_{kn}$ has such asymptotic starting at $n^{-k}$, as already remarked right after we obtained (\ref{eq:Legendre_asymptotics}).)

But if $n$ large enough, depending only on $m$, so that $t \le \log^2n$ also makes $t/2n$ small enough, we have uniformly in $t \in [0, \log^2 n]$
$$
(2n+1)\log \left( 1-(1-e^{-t/(2n)}) \tfrac{d}{d+1} \right) =
    - \frac{d}{d+1}t +\frac{p_1(t)}{n}+ \cdots + \frac{p_m(t)}{n^m}+O_{d,m}(n^{-m-1}),
$$
with fixed polynomials $p_k(t)$, where the $d$ dependence in each is of the type $(\frac{d}{d+1})^r$ for various $r$ and we notice that in each term $(\frac{d}{d+1})^r t^{r_1}$ we actually have $r \le r_1$.
For example $4p_1(t)=-(\frac{d}{d+1})^2 t^{2}+\frac{d}{d+1} t^{2}-2\frac{d}{d+1} t$.

Hence the central approximation of the main term of each integrand in (\ref{eq:Idjk-tilde}) is of the form
$$
e^{-dt/(d+1)}\left( 1+\frac{\tilde p_1(t)}{n}+ \cdots +
    \frac{\tilde{p}_m(t)}{n^m} + O_{d,m}(n^{-m-1}) \right),
$$
where $\tilde{p}_k(t)$ has the same properties as $p_k(t)$ and clearly $p_1=\tilde{p}_1$.

On the other hand the other exponential terms in (\ref{eq:Idjk-tilde}),
$(e^{-t/(2n)}d-1)^{2k} \cdot e^{t/(4n)+tk/(2n)}$, have developments of the type $(d-1)^{2k}+s_1(t)/n+\dots+s_m(t)/n^{m}+O_{d,m}(n^{-m-1})$, but now the dependence in $d$ of the polynomials $s_k(t)$ is also polynomial only, so overall the integrals in (\ref{eq:Idjk-tilde}) are a sum of integrals of the type
$e^{-dt/(d+1)} \cdot t^{(d-1)/2+j-1}t^r$
with coefficients depending on $d$ as discussed above.

\vspace{.3cm}
Now when we complete the tails, since the main factor is $e^{-dt/(d+1)}$ we need to make the change of variable $u=\frac{d}{d+1} t$, which means that every power of $t$ in the above integrands generates a factor of $(\frac{d+1}{d})^{(d-1)/2+j+r}$, as well as a $\Gamma ((d-1)/2+j+r)$. But because the respective outer coefficient contains the term $1/\Gamma ((d-1)/2+j)$, the quotient of the respective two Gammas is a polynomial in $d$.

If we look at the main asymptotic term, when $k=j=0$, there is only one corresponding integrand $e^{-dt/(d+1)} \cdot t^{(d-1)/2-1}$. Keeping in mind that $q_n \sim \frac{1}{\sqrt {\pi n}}$, as well as that the asymptotic requires only $(d+1)^{2n}$, so we have an extra $d+1$ factor from (\ref{eq:Idjk-tilde}), and that the Gamma factor cancels out completely, we obtain $a_{d+1}= a_d \frac{(d+1)^{(d+1)/2}}{2 \, d^{d/2}}$. Given that $a_2=1$, we immediately obtain the formula claimed in Theorem~\ref{T:TheoremA}(i):
$$ a_d=\frac{d^{d/2}}{2^{d-1}}. $$
Putting all of these observations together, Theorem~\ref{T:TheoremA}(i) follows. Explicit computations of $r_m(d)$, for $m=1,2,3,4$ are available upon request.

\vspace{.3cm}
We note that for any $m$, the only corresponding asymptotic term that appears in (\ref{eq:xns-expanded3}) for $j=m$ is when $k=0$, so as before we have only one corresponding integrand, $e^{-dt/(d+1)}t^{(d-1)/2-1+m}$, which gives an extra $(\frac{d+1}{d})^m$, so giving the term $r_m(d)(\frac{d+1}{d})^m$ in the recurrence for $r_m(d+1)$ that follows from (\ref{eq:Idjk-tilde}).
More generally, the discussion above about the shape of the coefficients of the powers of $t$ involved in each asymptotic term, especially the fact that any coefficient $d/(d+1)$ appears to a power at most the power of the corresponding $t$, so in the final result is canceled by the $(d+1)/d$ power coming from completing the tails, gives the recurrence:

\begingroup\makeatletter\def\f@size{11}\check@mathfonts
\begin{equation} \label{eq:r-recurrence}
r_m(d+1) = (\tfrac{d+1}{d})^m r_m(d)+ (\tfrac{d+1}{d})^{m-1} Q_{1,m}(d)r_{m-1}(d)+
    \dots +\tfrac{d+1}{d} Q_{m-1,m}(d)r_1(d)+Q_{m,m}(d),
\end{equation}
\endgroup
where $Q_{k,m}$ are rational functions in $d$ with denominator a power of $d$ only. A more careful look at all the expansions involved in (\ref{eq:xns-expanded3}) actually shows that the degree of the denominators of $Q_{k,m}(d)$ are at most $k$, since they all come only from the terms of degree $k$ in $1/n$ in the expansion of (\ref{eq:Idjk-tilde}), and the Legendre recurrence terms are actually the only ones giving denominators powers of $d$, highest possible being $k$ for the $1/n^k$ error term there. This is because the carryover term from $r_{m-k}(d)$ corresponding to the extra $t^{m-k}$ in the starting integral $(15)$ is precisely reflected in the coefficient $(\frac{d+1}{d})^k$, while simplification of the corresponding Gamma factors gives only polynomials in $d$ as noted above. Also all $Q_{k,m}(d)$ have rational coefficients since $q_n={2n \choose n}/2^{2n}$ and 
$q_{jn}$ do when expanded in powers of $1/n$ (as of course have the Taylor series of the elementary functions involved in (\ref{eq:Idjk-tilde})).

For example $Q_{1,1}(d)=-\frac{1}{8d}$ so $r_1(d+1)=\frac{d+1}{d} r_1(d)-\frac{1}{8d}$, with $r_1(1)=0$, $r_1(2)=-1/8$. This results in $r_1(d)=\frac{1-d}{8}$.
Similarly $Q_{1,2}(d)=\frac{3}{8d}$ and $Q_{2,2}(d)=\frac{d^2}{192} + \frac{d}{96} + \frac{11}{192} - \frac{1}{96d} - \frac{7}{128d^2}$, giving $r_2(d+1)=(\frac{d+1}{d})^2r_2(d)+\frac{d+1}{d}\frac{3}{8d}\frac{1-d}{8}+\frac{d^2}{192} + \frac{d}{96} + \frac{11}{192} - \frac{1}{96d} - \frac{7}{128d^2}$, with $r_2(1)=0$, $r_2(2)=\frac{1}{128}$. This results in $r_2(d)=\frac{(d^2-1)(2d-3)}{384}$, as claimed.

\vspace{.3cm}
Finally, since
$\displaystyle{
A_{2n}^{(d)} = {\small {2n \choose n}} x_n^{(d)} =
  2^{2n} \, q_n \, x_n^{(d)}
}$,
part (ii) of Theorem~\ref{T:TheoremA} follows.

\subsection{Asymptotics of the first returns $B_{2n}^{(d)}$}\label{subsec:Bs}

Now we continue with the deeper discussion of the asymptotic for $B_{2n}^{(d)}$'s and related matters, which involve the singularity analysis method from \cite[ch.VI]{FS} as well as the theory of holonomic ($D$-finite) functions and $P$-recursive sequences (\cite{Stanley1,Stanley2}), with forays in the theory of Fuchsian ODE's (\cite{Wasow}, \cite{Ince}) and in particular Heun equations (\cite{Ronveaux}). In the process, we will show that the $A$ and $B$ sequences are fundamentally different from a combinatorial viewpoint, as the first are $P$-recursive but the second are not, except in dimension $1$. We will first state the main theorems and then we will indicate how to prove them.

\begin{theorem}\label{T:TheoremB}
For $d \ge 3$, let $m_d=\sum_{n \ge 0} \frac{A_{2n}^{(d)}}{(2d)^{2n}}< \infty$ denote the expected number of returns to the origin, plus one, while for $d \ge 5$, let $\tilde m_d=\sum_{n \ge 0} \frac{nA_{2n}^{(d)}}{(2d)^{2n}}< \infty$.

\begin{enumerate}
\item[(i)] For $d \ge 3$ odd, we have

\begin{equation} \label{eq:Bs-odd-asyms} 
B_{2n}^{(d)} = b_d\frac{(2d)^{2n}}{(\pi n)^{d/2}} \left(1+ \frac{b_1(d)}{n}+\frac{b_2(d)}{n^2}+\dots +\frac{b_m(d)}{n^m}+
O_{d,m}\left( \frac{1}{n^{m+1}} \right) \right),
\end{equation}
where $\displaystyle{ b_d=\frac{a_d}{m_d^2}= \frac{d^{d/2}}{m_d^2 \cdot 2^{d-1}}}$, for all $d \ge 3$.
Here $b_1(3)=-\frac{3}{16}+\frac{9}{32m_3}-\frac{81}{16\pi^2m_3^3}$ and $b_1(d)=-\frac{d}{8}-\frac{d \, \tilde m_d}{m_d}$, for $d \ge 5$.

\item[(ii)] For $d \ge 4$ even,  we have

\begingroup\makeatletter\def\f@size{10}\check@mathfonts
\begin{equation} \label{eq:Bs-even-asyms}
B_{2n}^{(d)} = b_d\frac{(2d)^{2n}}{(\pi n)^{d/2}} \left(1+ \frac{b_1(d, \log n)}{n}+\frac{b_2(d, \log n)}{n^2}+\dots +\frac{b_m(d, \log n)}{n^m}+
O_{d,m}\left( \frac{b_{m+1}(d, \log n)}{n^{m+1}} \right) \right)
\end{equation}
\endgroup
where $b_k(d, \log n)$ are now (possibly constant) polynomials in $\log n$, with the first non-constant polynomial $b_m(d, \log n)$ for $m=d/2-1$ with leading term $-2{d-2 \choose d/2-1}\frac{a_d}{\pi^{d/2}m_d} \log n$.

In particular, for $d=4$ we have
\begin{equation} \label{eq:B4-asyms}
B_{2n}^{(4)}=b_4 \frac{8^{2n}}{(\pi n)^{2}}
    \left( 1-\frac{8}{\pi^2 m_4}\frac{\log n}{n} + O\left( \tfrac{1}{n} \right) \right),
\end{equation}
while for $ d \ge 6$ we have $b_1(d, \log n)=b_1(d)=-\frac{d}{8}+\frac{d \, \tilde m_d}{m_d}$.

\item[(iii)] For $d=2$, we have (\cite[p.426]{FS}, with $\gamma$ the Euler's constant and $K\doteq 0.8825424$)
\begin{equation} \label{eq:B2-asyms}
B_{2n}^{(2)}=\pi \frac{4^{2n}}{n \log^2 n}
    \left( 1 - 2 \frac{\gamma + \pi K}{\log n} + O\left( \frac{1}{\log^2 n} \right) \right).
\end{equation}

\item[(iv)] For $d=1$, we have
\begin{equation} \label{eq:B1-asyms}
B_{2n}^{(1)} = \frac{A_{2n}^{(1)}}{2n-1} = \frac{2^{2n}}{2n \sqrt {\pi n}}
    \left( 1-\frac{3}{8n}+\frac{25}{128 n^2} + O\left( \frac{1}{n^3} \right) \right).
\end{equation}

\end{enumerate}
\end{theorem}

After recalling a few definitions and making some observations, we collect all the information about the generating functions of $x_n^{(d)}$, $A_{2n}^{(d)}$, and $B_{2n}^{(d)}$ in Theorems~\ref{T:TheoremXs}, \ref{T:TheoremAs}, and \ref{T:TheoremBs}, respectively. Their proofs will be given in the next subsection, but here we show how they imply Theorems \ref{T:TheoremC} and \ref{T:TheoremB}, using singularity analysis theory results like the zig-zag algorithm from \cite[Ch.VI]{FS}.

\vspace{.3cm}
If $f=\sum a_n w^n$ and $g=\sum b_n w^n$ are Taylor series with (nonzero) radii of convergence $R_1$ and $R_2$, respectively, we define their {\it Hadamard convolution} by
$h(w)=f*g(w)=\sum a_n b_n \, w^n$.
From the integral representation
$$
h(w)=\frac{1}{2\pi i}\int_{|\zeta|=r}f(\zeta)g(w/\zeta)\, \frac{d\zeta}{\zeta},
    \text{ for } |w| <R_1R_2, |w|/R_2 <r<R_1,
$$
$h$ has radius of convergence at least $R_1 R_2$.

\vspace{.3cm}
Following \cite[App.B.4]{FS}, a formal power series over some field $K$ (for us $K=\mathbb C$, $K=\mathbb Q$ or $K=\mathbb Z/p$, $p$ prime) $f(z)=\sum_{n \ge 1} f_n z^n$ is called {\it holonomic} or {\it $D$-finite} (short for differentiably finite) if the vector space over $K(z)$ (the field of rational functions with coefficients in $K$) spanned by the set of all its formal derivatives $(\partial^j f(z))_{j=0}^{\infty}$ is finite dimensional.
This is equivalent to the existence of a (non-trivial) finite degree linear ODE with coefficients in $K(z)$ (hence by clearing denominators in $K[z]$) formally satisfied by $f$.

Similarly, a sequence $(f_n)_n \in K$ will be called {\it $P$-recursive} (short for polynomially recursive) if there are polynomials $P_0, \dots, P_r \in K[z]$ such that
$\sum_{k=0}^{r} P_k(n) \, f_{n+k}=0$, for $n \ge n_0$.
It is a fairly simple matter to note that the holonomicity of a generating formal power series $f$ is equivalent to the $P$-recursiveness of its coefficients $f_n$.

If $f$ is algebraic over $K$, so there are polynomials $Q_0, \dots , Q_r \in K[z]$, not all zero, such that $\sum_{k=0}^r Q_k(z)f^k(z)=0$, then $f$ is holonomic (all the formal derivatives of $f$ belong to $K(z,f)$ which is finitely dimensional over $K(z)$ since $f$ is algebraic). If $f$ is not algebraic, we call it transcendental. Rational functions like $\frac{1}{1-z}=\sum_{n \ge 0} z^n$ are algebraic over any field $K$, while $e^z=\sum z^n/n!$ is transcendental over $\mathbb C$ but holonomic since $(e^z)'-e^z=0$, or, equivalently, $(n+1)\cdot \frac{1}{(n+1)!}-\frac{1}{n!}=0$ being the $P$-recurrence for the coefficients. By extension, we will call the sequence of coefficients $(f_n)_n$ algebraic or transcendental over $K$ if the corresponding formal series is, though of course the definition is useful only for $P$-recursive coefficients as non-$P$-recursive coefficients cannot be algebraic anyway. In particular if the $f_n$'s are integers, it makes sense to talk about algebraicity modulo a prime $p$ in addition to algebraicity over $\mathbb C$. By Eisenstein's Lemma (\cite[part VIII, no.150]{PS}), if $f$ is an algebraic power series over $\mathbb C$ with rational coefficients, so $P_0(z)+P_1(z)f(z)+...+P_n(z)f^n(z) =0$, $n \ge 1$, $P_n \ne 0$, $P_k \in \mathbb C[z]$, then we can actually take $P_k \in \mathbb Z[z]$.

The fundamental difference between $\mathbb C$ and $\mathbb Z/p$ for integral sequences is that pointwise multiplication of coefficient sequences $(f_n)\, (g_n) \to (f_n g_n)$\footnote{Over the complex field, this pointwise multiplication corresponds to Hadamard convolution, if the formal power series have non  zero radius of convergence, and by a slight abuse of terminology is also called convolution or even Hadamard convolution over any field $K$, though we do not have an integral representation in general.} preserves algebraicity only in finite characteristic $p$, though it preserves holonomicity for any field (see \cite[Ch.6]{Stanley2}).

\vspace{.3cm}
If we work over $\mathbb C$ and $f=\sum f_nz^n$ holonomic has non-zero radius of convergence, so it actually defines an analytic function and not only a formal power series, we can use all the machinery of linear ODEs with polynomial (rational) coefficients (\cite{Ince}, \cite{Wasow}). In particular such an ODE is called {\it Fuchsian} if all its solutions are locally bounded by rational functions. Equivalently, the ODE is Fuchsian iff when denoting the singularities of the ODE (the poles of the coefficients $R_k$ when the equation is written as
$f^{(m)}(z)+\sum_{0 \le k \le m-1} R_k(z) \, f^{(k)}(z)=0$,
with $R_k$ rational functions) by $q_1, \dots, q_M$ (including $\infty$ if when the equation is transformed by $ z \to 1/z$, we get a singularity at $0$), then, for all local solutions of the ODE, there is a $N$ such that $(z-q_k)^N \, f$ is bounded near $q_k$, with the corresponding definition at $\infty$. The usual definition (\cite{Ince}) involves the notion of regular singular point but is equivalent to the above. The Fuchsian condition constrains strongly the shape of the equation (\cite{Forsyth}). Also a linear ODE with apriori meromorphic coefficients only, which is Fuchsian in the sense above, can actually be taken to have rational coefficients, while if $f$ is a solution of a linear ODE (again apriori with meromorphic coefficients only) and all the singularities of $f$ are Fuchsian (hence they are finitely many), then $f$ is is the solution of a Fuchsian ODE of degree at most the degree of the original ODE satisfied by $f$. In particular, algebraic functions over $\mathbb C$ (which are automatically locally analytic except at their finitely many singularities) satisfy Fuchsian ODEs, though it is a classical and quite difficult problem to determine if a given Fuchsian ODE has at least one algebraic solution (\cite{Forsyth}).

\vspace{.3cm}
Conform \cite{FS} or \cite{FFK}, an analytic function $f(w)$ given by a Taylor series with radius of convergence $1$ is called {\it $\Delta$-analytic} (or {\it $\Delta$-regular}) if it has an analytic extension to a $\Delta$-domain of the type $|w|<1+\eta$, $| \arg (w-1) | > \phi$, for some $\eta>0$ and $0< \phi <\pi/2$. More generally we call any analytic $f$ given by a Taylor series with arbitrary finite radius of convergence $\Delta$-analytic if $f(\zeta w)$ is so in the restricted sense above for some complex $\zeta \ne 0$. The theory of linear ODEs with polynomial (rational) coefficients (\cite{Ince}) shows that holonomic functions defined by convergent power series and being solutions of such ODEs are $\Delta$-analytic, though of course the converse if far from true since $\Delta$-analyticity depends only on the singularities on the boundary of the disc of convergence.

\vspace{.3cm}
We define $T(w)=1-w$ and $L(w)=-\log (1-w)=\sum_{n \ge 1}\frac{w^n}{n}$, for $|w|<1$.

\vspace{.3cm}
$\Delta$-analyticity is a very useful property enjoyed by many elementary functions like $T^a=(1-w)^a$ (algebraic), $(L/w)^b=(\log (1/(1-w))/w)^b$, for $a,b \in \mathbb R$, and polylogarithms $\sum_{n \ge 1}\frac{w^n}{n^k}$, $k \ge 2$, while as noted holonomy is much more restrictive. For example, while polylogarithms and $T^a$ are holonomic, only $L^n$, $n \ge 1$ are holonomic. While both properties are preserved under various operations (\cite[VI]{FS}, \cite{Stanley2}), including differentiation, integration, Hadamard convolution $f*g$, and multiplication $f \, g$, $\Delta$-analyticity is preserved under division $f/g$, with $g$ having no zeroes on a $\Delta$-domain, but holonomy in general is not preserved under division. For example, if $f$ is holonomic ($f(0) \ne 0$), the theorem of Harris and Sibuya \cite{HS} states that $1/f$ is holonomic if and only if $f'/f$ is algebraic.

\vspace{.3cm}
Our next theorems state that the generating functions of the first returns are also $\Delta$-analytic, but only the $A$ and $x$ sequences lead to holonomic functions and not the $B$ sequence.

\begin{theorem}\label{T:TheoremXs}
For
$\displaystyle{
F_d(w)=\sum_{n \ge 0} x_n^{(d)} \, w^n,
}$
the following hold:

\begin{enumerate}
\item[(i)] $F_d$ has radius of convergence $1/d^2$ and is uniformly convergent on the circle of convergence for $d \ge 4$. It is holonomic, hence $\Delta$-analytic, for all $d \ge 1$ and satisfies a Fuchsian linear ODE (which we conjecture to be of degree $d-1$). In particular, $F_d$ can be extended to a multivalued analytic function in the plane with finitely many singularities.

\item[(ii)] For $d \ge 2$, $F_{d}$ has finite order singularities at $1/k^2$, for all $1 \le k \le d$ and $k$ having the same parity as $d$, and at $\infty$, with principal branch extending the power series on the plane cut from $1/d^2$ to $\infty$.

\item[(iii)] All the finite singularities of $F_d$ are of the same type and with the same kind of asymptotics (which is also a development on angular domains) as described below in Theorem 7 for the circle of convergence singularity. At infinity the singularity will always be of a logarithmic nature (for $d\geq 3$).

\item[(iv)] $\displaystyle{ F_1(w)=1/(1-w)}$, $\displaystyle{ F_2(w)=1/\sqrt {1-4w} =A_1(w)}$, $F_3(1/9)=\infty$.

\end{enumerate}
\end{theorem}

\begin{theorem}\label{T:TheoremAs}
For
$\displaystyle{
A_d(w)=\sum_{n \ge 0} A_{2n}^{(d)} \, w^n,
}$
the following hold:

\begin{enumerate}
\item[(i)] $A_d(w)=F_d(w)*A_1(w)=F_d(w)*F_2(w)$, for all $d \ge 1$.
$A_1(w)=1/\sqrt {1-4w}=F_2(w)$ and $A_2(w)=A_1(w)*A_1(w)= \ \tensor[_2]{F}{_1} \left( \tfrac{1}{2}, \tfrac{1}{2}; 1; 16w \right)$ is a hypergeometric function.

\item[(ii)] $A_d$ has radius of convergence $1/(2d)^2$ and is uniformly convergent on the circle of convergence for $d \ge 3$, with $A_d(1/(2d)^2)=m_d >1$, for $d \ge 3$, and $A_1(1/4)=A_2(1/16)=\infty$. It is holonomic, hence $\Delta$-analytic, for all $d \ge 1$, and satisfies a Fuchsian linear ODE (which we conjecture to be of degree $d$). In particular, $A_d$ can be extended to a multivalued analytic function in the plane with finitely many singularities.

\item[(iii)] $A_d$ has finite order singularities at $1/(2k)^2$, for all $1 \le k \le d$ and $k$ having the same parity as $d$, and at $\infty$, with its principal branch extending the power series on the plane cut from $1/(2d)^2$ to $\infty$.

\item[(iv)] All the finite singularities of $A_d$ are of the same type and with the same kind of asymptotics (which is also a development on angular domains) as described below in Theorem 7 for the circle of convergence singularity. For $d$ even, the singularity at infinity will always have logarithmic and half integral powers (e.g., it is well known that $A_2(w) \sim (c_1+c_2 \log w)\sqrt {1/w}$ for $|w| \to \infty$ \cite{DLMF}). For $d \ge 3$ odd, the singularity at infinity is of square root type and contains no logarithmic terms.

\end{enumerate}
\end{theorem}

\begin{theorem}\label{T:TheoremBs}
For
$\displaystyle{
B_d(w)=\sum_{n \ge 1} B_{2n}^{(d)} \, w^n,
}$
the following hold:

\begin{enumerate}
\item[(i)]$(1-B_d(w)) \cdot A_d(w)=1$.

\item[(ii)] $B_d$ has radius of convergence $1/(2d)^2$ and is uniformly convergent on the circle of convergence. It is not holonomic for $d \ge 2$, but it is still $\Delta$-analytic, for all $d \ge 1$.

\item[(iii)] $B_1(w)=1-\sqrt {1-4w}$, $B_1(1/4)=B_2(1/16)=1$, and, for $d \geq 3$, $B_d(1/(2d)^2)=1-1/m_d=p_d$, P\'{o}lya's return probability.

\end{enumerate}
\end{theorem}

\vspace{.3cm}
Next we consider the normalized series
$\tilde{F_d}(w)=\sum_{n \ge 0}\frac{x_n^{(d)}}{d^{2n}} \, w^n$,
with radius of convergence~1, and the analogue definitions for $\tilde{A_d}(w)$ and $\tilde{B_d}(w)$.

\begin{theorem}\label{T:TheoremC}
The following hold:
\begin{enumerate}
\item[(i)] For $d \ge 2$ even, $\tilde {F_d}$ has an asymptotic expansion in terms of half integer and integer powers of $T$, with the first asymptotic term being
    $\frac{a_d \, \Gamma(\frac{3-d}{2})}{\pi^{(d-1)/2}} T^{(d-3)/2}$.
    For $d \ge 5$ odd, $\tilde {F_d}$ has an asymptotic expansion in terms of $LT^k$ and $T^m$, with the first asymptotic term being $\frac{a_d}{d! \, \pi^{(d-1)/2}} LT^{(d-3)/2}$. Also
    $\tilde F_1=1/T$ and $\tilde{F_3}=\frac{3\sqrt 3}{4\pi}L + f_{3,1}+ f_{3,2}LT+\cdots$.

\item[(ii)] For $d \ge 1$ odd, $\tilde{A_d}$ has an asymptotic expansion in terms of half integer and integer powers of $T$, with the first asymptotic term being
    $\frac{a_d \, \Gamma(\frac{2-d}{2})}{\pi^{d/2}} T^{(d-2)/2}$.
     For $d \ge 4$ even, $\tilde{A_d}$ has an asymptotic expansion in terms of $LT^k$ and $T^m$,  with the first first asymptotic term being $\frac{a_d}{d! \, \pi^{d/2}} LT^{(d-2)/2}$. Also $\tilde{A_2}=L/\pi + K+ a_{2,2}LT+\cdots$.

\item[(iii)] For $d \ge 1$ odd, $\tilde{B_d}$ has an asymptotic expansion in terms of half integer and integer powers of $T$, with the first asymptotic term being
    $\frac{a_d \, \Gamma(\frac{2-d}{2})}{m_d^2 \, \pi^{d/2}} T^{(d-2)/2}$.
     For $d \ge 4$ even, $\tilde{B_d}$ has an asymptotic expansion in terms of $L^j T^k$, $k-j \ge d/2-2$, and $T^m$, $k,j,m \ge 1$, with the first first asymptotic term being
     $\frac{a_d}{d! \, \pi^{d/2} \, m_d^2} LT^{(d-2)/2}$.
     For $d=2$ the asymptotic is in terms of $L^{-j}$, $j \ge 1$, since they are all higher order (converging much slower to zero at $1/16$) than any $L^j T^k$, $j \ge 0$, $k \ge 1$.

\end{enumerate}
\end{theorem}

The fundamental result is Theorem~\ref{T:TheoremXs} and most notably the statements about holonomicity and the principal branch which we will prove rigorously in the next subsection, while we will give a sketch of proof for the rest of the claims. Theorems \ref{T:TheoremAs} and \ref{T:TheoremBs} follow from Theorem~\ref{T:TheoremXs} and a few basic considerations which we will also address in the next subsection. Here we prove Theorems~\ref{T:TheoremC} and \ref{T:TheoremB}, given Theorems \ref{T:TheoremXs}, \ref{T:TheoremAs}, and \ref{T:TheoremBs}.

We have explicit computations in dimension $1$, since we know the $A$ series exactly to be $1/\sqrt {1-4w}$, so we can compute the $B$ series and hence its coefficients from $(1-B) \cdot A=1$. When $d=2$, the random $2$-walk is the direct product of two independent one dimensional random walks and the full analysis of this fairly difficult case, as the small logarithmic error terms show, is done in \cite[Example VI.14]{FS}. These remarks justify parts (iii) and (iv) of Theorem~\ref{T:TheoremB}, which we included in the statement for completion, so we prove Theorems~\ref{T:TheoremB} and \ref{T:TheoremC} only in the case $d \ge 3$.

Note that for $d \ge 3$ the random walk we consider is not a product of $1$-dimensional walks and we see this reflected in the formulas for $A_{2n}^{(d)}$'s and $B_{2n}^{(d)}$'s, which involve $(2d)^{2n}$ rather than $2^{2nd}$. However, due to the similar asymptotic expansion of the Legendre polynomials $(d^2-1)^nP_n(\frac{d^2+1}{d^2-1})$ to the one for the binomial coefficient $2n \choose n$, the asymptotic expansions of $A_{2n}^{(d)}/(2d)^{2n}$ and $B_{2n}^{(d)}/(2d)^{2n}$ are similar to the ones from the $d$-product of $1$-dimensional walks discussed by Flajolet and Sedgewick (P\'{o}lya's drunkard problem referenced above).

Theorems~\ref{T:TheoremXs} and \ref{T:TheoremAs} (and their proofs) show inductively that $F$ and $A$ (hence $\tilde F$ and $\tilde A$) are analytic functions amenable to singularity analysis in the sense of \cite[Ch.VI]{FS}, so they are not only holonomic, hence $\Delta$-analytic, but also admit an asymptotic representation in terms of the basic functions $T, L$ and their powers which can be determined (up to non-negative powers of $T$) from the asymptotics of their coefficients.

The coefficients of $\tilde{A_d}(w)$ are asymptotic with
$$
\sum_{k \ge 1} a_{d,k}\frac{1}{(\pi n)^{d/2+k-1}}
$$
and those of $\tilde {F_d}(w)$ are asymptotic with
$$
\sum_{k \ge 1} f_{d,k}\frac{1}{(\pi n)^{(d-1)/2+k-1}},
$$
where $f_{d,1}=a_{d,1}=\frac{d^{(d/2)}}{2^{d-1}}$. We know that $\tilde{F_1}(w)=(1-w)^{-1}=T^{-1}$ and $\tilde {F_2}(w)=\tilde {A_1}(w)=(1-w)^{-1/2}$, hence when the coefficients are asymptotic with odd powers of $1/n$ the function is asymptotic to a sum of half integer (and integer) powers of $T$, and when the dominant term is $n^{q/2}$, $q \ge 3$, the series being convergent, we start with the constant $\tilde{A_d}(1)=m_d$, $d \ge 3$ odd, and $\tilde{F_d}(1)=f_d$, $d \ge 4$ even.

When the powers are even and the dominant term is $1/n$ the series are divergent at $1$ and start with $L/\pi$ for $\tilde{A_2}(w)$, $a_2=1$, and $\frac{3\sqrt 3}{4\pi}L$ for $\tilde {F_3}(w)$, $a_3=f_3=\frac{3\sqrt 3}{4}$, respectively, followed by a sum of $T^j$, $j \ge 0$, and $LT^m$, $m \ge 1$. When the dominant term is $1/n^q$, $q \ge 2$, the series are again convergent at $1$. This means that they start with a constant $m_d$ or $f_d$ and then we have a sum of $T^j$, $j \ge 1$, and $LT^m$, $m \ge q-1 \ge 1$. These comments prove the claims from Theorem~\ref{T:TheoremC} about the first ``true" asymptotic term being $LT^{(d-2)/2}$ for $\tilde{A_d}(w)$ and $LT^{(d-3)/2}$ for $\tilde{F_d}(w)$, respectively.

Since the asymptotic of $\tilde{B_d}(w)$ follows from the relation
$
 (1-\tilde{B_d}(w)) \cdot \tilde{A_d}(w)=1
$,
essentially by identifying coefficients, we have to do an analysis considering how fast terms converge to $0$ at $1$, so we have to look at the terms $O(1)$, $o(T)$, $O(T)$, $o(T^2)$, $O(T^2)$, and so on. There are two cases to be considered.

\vspace{.3cm}
\noindent
{\it Case 1.} When $d \ge 3$ odd (using that $T^{d/2}$ has main asymptotic coefficient $\frac{1}{\Gamma(\frac{2-d}{2})}$),
$$
  \tilde{A_d}(w)=m_d+\sum_{k \ge 1}\tilde{a}_{d,k} T^{d/2+k-2}+\sum_{j \ge 1}q_{d,j}T^j, \text{ with }
    \tilde{a}_{d,1} = \frac{a_d \, \Gamma(\frac{2-d}{2})}{\pi^{d/2}},
$$
hence by identifying coefficients we get
$$
  \tilde{B_d}(w)=(1-1/m_d) + \sum_{k \ge 1} \tilde{b}_{d,k} T^{d/2+k-2}+\sum_{j \ge 1}\tilde{q}_{d,j}T^j, \text{ with }
    \tilde{b}_{d,1} = \frac{a_d \, \Gamma(\frac{2-d}{2})}{\pi^{d/2} \, m_d^2}.
$$
In particular the dominant power of the coefficients of $\tilde {B_d}(w)$ has coefficient $ \frac{a_d}{(\pi n)^{d/2} \, m_d^2}$, which also gives the main claim of Theorem~\ref{T:TheoremB} for $d \ge 3$ odd, as the error terms are the same as for the $A$ coefficients.

\vspace{.3cm}
\noindent
{\it Case 2.} When $d \ge 4$ even (using that $LT^{d/2}$ has main asymptotic coefficient $d!$),
$$
  \tilde{A_d}(w)=m_d+\sum_{k \ge 1}\tilde{a}_{d,k} LT^{d/2+k-2}+\sum_{j \ge 1}q_{d,j}T^j,
   \text{ with } \tilde{a}_{d,1}=\frac{a_d}{\pi^{d/2} \, d!}.
$$
Here we notice that, when we write $(1-\tilde {B_d}(w))\cdot \tilde {A_d}(w)=1$, the constant and the first true asymptotic term are the same\footnote{We also have powers of $T$ of order less than $LT^{d/2-2}$ for $d \ge 6$ say when we start with $LT^{m}$, $m \ge 2$, so have $T$, \dots , $T^{m-1}$ before that in asymptotic order at $1$, but those are of course negligible as the coefficients asymptotic go since they change only the first $m-1$ coefficients.} as in the odd case, so

$$
 \tilde{B_d}(w)=(1-1/m_d)+\sum_{j \ge 1}\tilde{q}_{d,j}T^j+\frac{1}{\pi^{d/2} \, d!} \frac{a_d}{m_d^2} LT^{d/2-1}+\cdots,
$$
But now we notice that when we do the term by term multiplication, we start having also terms $L^2T^j$, $L^3T^j$, and so on, first coming from the need to cancel the $LT^j$ terms in $A$ multiplied with a $LT^k$ term in $B$ and then of course we need to cancel $LT^j$ terms in $A$ multiplied with the previous $L^2T^k$ in $B$ etc. Conform \cite[Ch.VI]{FS} we know that these terms add logarithmic errors with the logarithmic term a polynomial of degree $m-1$ and the dominant power of $1/n$ being as usual $k+1$ in $L^mT^k$, so in particular for $L^2T^j$ we get a dominant term $\log n /n^{j+1}$. But now since the first asymptotic term is $LT^{d/2-1}$ the first $L^2T^j$ term is $L^2T^{d-2}$ with dominant coefficient term $\log n/{n^{d-1}}$. On the other hand the second asymptotic term in $A$ is $LT^{d/2}$, hence $B$ also has this term $LT^{d/2}$ with ``good" dominant term of order $1/n^{d/2+1}$. It follows that the only case when the $L^2T^{d-2}$ dominates and appears as main error is when $d=4$, so $d-2=d/2$. In this case $L^2T^2$ goes slower to zero at $1$ than $LT^2$, so it is of lower order and appears first in the asymptotic expansion and its main term $\log n/n^3$ becomes the main error here. For $d \ge 6$, so satisfying $d/2<d-2$, the term $LT^{d/2}$ remains the second asymptotic term, so its dominant term $1/n^{d/2+1}$ is still the main error in the asymptotic for the coefficients $B$ and with this Theorem~\ref{T:TheoremC} is proved.

Theorem~\ref{T:TheoremB} just restates the functional asymptotics of $\tilde B$ from Theorem~\ref{T:TheoremC} into asymptotics for the coefficients so the shape of the coefficients is the same. The explicit computations follow from simple algebraic manipulations using our earlier result that $a_1(d)=-d/8$, noting that $L^2T^k$ has main asymptotic coefficient $2k! \log n/n^{k+1}$, and the fact that for $d=3$ we have the explicit computation from \cite{Joyce72}, where the coefficients of $T^{1/2}$, $T^{3/2}$ follow from our considerations above too, but the coefficient of $T$ follows only from the method of the cited paper:

$$\tilde{A_3} = m_3 -\frac{3\sqrt 3}{2\pi}T^{1/2}+\frac{9}{32}(m_3+\frac{6}{\pi^2m_3})T-\frac{3\sqrt 3}{4\pi}T^{3/2} +O(T^2).$$

\noindent With this the proofs of Theorem~\ref{T:TheoremB} and Theorem~\ref{T:TheoremC} are completed.

\subsection{Analytic Proofs }\label{subsec:AnP}

\begin{lemma}\label{L:Lemma TT}
The following identities hold for $|w|<1$:
$$
\sum_{m=0}^{\infty} {n+m \choose m}^2 w^m=(1-w)^{-2n-1}\sum_{k=0}^{n} {n \choose k}^2 w^k = (1-w)^{-1}(w-1)^{-n} P_n\left( \frac{w+1}{w-1} \right).
$$
\end{lemma}

\begin{proof}
The first equality can be found in \cite{Turan,Takacs}. The last equality follows from our earlier formula (\ref{eq:Legendre_recurrence}) and the fact that $P_n(-x)=(-1)^n P_n(x)$.
\end{proof}

Hence when we substitute the recurrence $x_n^{(d+1)}=\sum_{k=0}^n{n \choose k}^2 x_k^{(d)}$ in the series for $F_{d+1}$, we obtain from Lemma~\ref{L:Lemma TT} that:

\begin{equation} \label{eq:Conv0}
F_{d+1}(w)=\sum_{n \ge 0} \left( x_n^{(d)}(1-w)^{-1}(w-1)^{-n} P_n\left( \frac{w+1}{w-1} \right)\right) w^n,
\text{ for }|w|< 1/(d+1)^2.
\end{equation}

\noindent
Denoting $u=\frac{w}{w-1}$, such that $(w-1)(u-1)=1$ and $\frac{w+1}{w-1}=2u-1$, the equation \eqref{eq:Conv0} becomes:

\begin{equation} \label{eq:Conv1}
F_{d+1}(w)=(1-u)\, \sum_{n \ge 0} (x_n^{(d)}\, P_n(2u-1)\, u^n),
\text{ for }|w|< 1/(d+1)^2.
\end{equation}

By the Legendre polynomials generating function definition, we recall that
\begin{equation} \label{eq:Leg1}
G(t,u)=\sum_{n \ge 0}P_n(u)t^n=(1-2tu+t^2)^{-1/2}.
\end{equation}

\noindent
Writing $P_n(2u-1)=\sum_{m \ge 0}q_{m,n} \, u^m$, we consider two power series in two variables (convergent on a small neighborhood of $(0,0)$):

\begin{align*}
V_1(t,u) &= \sum_{m,n \ge 0} x_n^{(d)} \, t^n u^m = \frac{1}{1-u} F_d(t), \\
V_2(t,u) &= \sum_{m,n \ge 0} q_{m,n} \, t^n u^m = \sum_{n \ge 0} P_n(2u-1)\, t^n = G(t,2u-1).
\end{align*}
With this notation, Equation \eqref{eq:Conv1} can then be written as

\begin{equation} \label{eq:ConvM}
(1-w)F_{d+1}(w)= (V_1*V_2)(t,u)|_{t=u}=(\frac{1}{1-u}F_d(t))*G(t,2u-1)|_{t=u}, \text{ for }w=\frac{u}{u-1}.
\end{equation}

Since the kernel $G$ is algebraic and the two variable convolution evaluation is at the linear relation $t=u$, we have at our disposal the powerful convolution theorems for holonomic functions (\cite{Stanley2}, \cite{Lipshitz}). In particular since $F_1(w)=\frac{1}{1-w}$ is rational, we get inductively (\cite[Thm.2.7, Prop.2.3]{Lipshitz}) that $F_d(w)$ is holonomic for all $d \ge 2$. Also since $\frac{1}{1-t}$ is the convolution unit, we get that
$F_2(w)=\frac{1}{1-w}\, G(\frac{w}{w-1}, 2\frac{w}{w-1}-1)=(1-4w)^{-1/2}$, as expected.
Similarly, since $A_1=F_2$ and $A_d(w)=F_d(w)*A_1(w)=F_d(w)*F_2(w)$, for $d \ge 1$, we get inductively that $A_d$ is holonomic for all $d \ge 1$.

Using again that $\frac{1}{1-u}$ is the convolution unit (now with respect to $u$), we can explicit the formula \eqref{eq:ConvM}, for $u$ sufficiently small, as:

\begin{equation} \label{eq:ConvExp}
(1-w)F_{d+1}(w) = \frac{1}{2\pi i}\int_{|\zeta|=s_d}
    F_d(\frac{1}{\zeta})\, G(\zeta u, 2u-1)\, \frac{d\zeta}{\zeta},
     \text{ for } s_d > d^2, w=\frac{u}{u-1}.
\end{equation}
Changing variables back to $w$ and $t=\frac{1}{\zeta}$ we can also write

\begin{equation} \label{eq:ConvExpN}
F_{d+1}(w) = \frac{1}{2\pi i}\int_{|t|=1/s_d}
    F_d(t)P(t,w)^{-1/2}dt, \text{ for } P(t,w)=w^2-wt(2w+2)+t^2(w-1)^2.
\end{equation}

Next we will use the method of Hadamard (\cite[Article 88]{Dienes}) to analyze inductively the singularities of $F_{d+1}$ in terms of those of $F_d$ and $P^{-1/2}$ by deforming the circle $|t|=1/s_d$ in such a way that the singularities of $F_d(t)$ stay outside all the deformation curves, while the singularities of $P^{-1/2}$ stay inside.

We first notice that $P(w,t)=P(t,w)$ and $P$ has roots

\begin{equation} \label{eq:roots}
t_1(w)=\frac{w}{(\sqrt w +1)^2} \text{ and }
t_2(w)=\frac{w}{(\sqrt w -1)^2},
\end{equation}
where changing the determination of the square root, just switches them. Also $t_1(0)=t_2(0)=0$, $t_1(\infty)=t_2(\infty)=1$, $t_1(1)=\frac{1}{4}$, $t_2(1)=\infty$, in the usual limit sense in the latter cases (i.e., as $w \to 1$ the roots $t_1(w) \to \frac{1}{4}$, $t_2(w) \to \infty$). It is the case that $t_{1,2}(w) \in [0, \infty]$ if and only if $w \in [0, \infty]$.

Given a Jordan curve $t \to \Gamma(t)$ in the (extended) $t$ plane and taking $w$ in the (extended) $w$ plane, the only condition for a continuous (hence analytic) determination of $t \to P^{-1/2}(t,w)$ is that the zeroes $t_1(w)$, $t_2(w)$ are not separated by $\Gamma$ (and of course are not on $\Gamma$ either). Then, by fixing an initial determination at some $w_0$, such a determination can be chosen jointly continuous on any $\Gamma \times U$, where $U \ni w_0$ is a connected open set in the $w$ plane with the property above. In particular, if we can also choose a continuous determination of $F_d(t)$ on $\Gamma$ (where we proceed inductively to define $F_d$ as a multivalued analytic function starting from the given power series around zero), the corresponding formula

\begin{equation} \label{eq:ConvExpNC}
F_{d+1, \Gamma}(w) = \frac{1}{2\pi i}\int_{\Gamma}
    F_d(t)P(t,w)^{-1/2}dt, \text{ for } P(t,w)=w^2-wt(2w+2)+t^2(w-1)^2
\end{equation}
defines an analytic function for $w$ in the open set for which $t_{1,2}$ are not separated by $\Gamma$. However, if $\Gamma$ doesn't separate at least a singularity of $F_d$ from the zeroes of $P$, a simple deformation argument shows that for such $w$ the integral is zero, hence for any $\Gamma$ the set of interest $U_{\Gamma}$ is defined to contain only those $w$ for which $\Gamma$ separates the singularities of $F_d$ from the zeroes $t_{1,2}(w)$. In general $U_{\Gamma}$ is an open but possibly disconnected set and $0 \in U_{\Gamma}$ precisely when $\Gamma$ separates $0$ from the singularities of $F_d$ (since $t_{1,2}(0)=0$). In particular this means that $0$ is inside $\Gamma$ by the inductive description of those singularities. Clearly any $w$ for which $t_{1,2}(w)$ is not a singularity of $F_d$ will belong to such a $U_{\Gamma}$ since $F_d$ has only finitely many singularities and we claim that for any $w \ne 1, \infty$ we can find $\Gamma$ so that $0 \in U_{\Gamma}$ and $w$ is in the connected component of $U_{\Gamma}$ containing $0$. This then gives us an analytic continuation of $F_{d+1}$ from the power series at zero (given also by \eqref{eq:ConvExpN}) to $w$, as we can choose the circle around $0$ in \eqref{eq:ConvExpN} to be small enough so is contained inside $\Gamma$, and then Cauchy theorem immediately implies that
$$
\int_{|t|=1/s_d}F_d(t)P(t,w)^{-1/2}dt=\int_{\Gamma}F_d(t)P(t,w)^{-1/2}dt,
$$
for $|w|$ small enough, so $F_{d+1}(w)=F_{d+1,\Gamma}(w)$ for $|w|$ small enough.

The claim about singularities follows (inductively) since for any $w_0 \notin [0, \infty]$ we have $t_{1,2}(w_0) \notin [0, \infty]$, so taking a small tubular neighborhood $T_{w_0}$ of the segment $[0,w_0]$ that stays away from the positive half axis except near $0$, we have that $\{ t_{1,2}(w) | w \in T_{w_0} \}$ is a compact set that stays away from the positive half axis except for a small neighborhood of zero, while the singularities of $F_d$ are inductively on $[0, \infty]$ and at a positive distance from zero, so we can separate the two sets by a Jordan curve $\Gamma$.

Let $r\in [0, \infty)$,  $r \ne 1$, be such that $t_1(r)$, $t_2(r)$ are not among the inductively assumed (finite) singularities of $F_d$, so not among $\frac{1}{k^2}$, for $1 \le k \le d$ with $k$ and $d$ of same parity. We consider $T_r$ a small tubular neighborhood of the curve formed by the segments $[0,ir\sin \epsilon], [ir\sin \epsilon, ir \sin \epsilon+r \cos \epsilon]$ and the arc $[re^{i\epsilon}, r]$. For $r<1$ we choose $\epsilon <0$ such that $\rho < \cos \epsilon$ for all $\rho>0$ in the tubular neighborhood, while for $r>1$ we can choose any $\epsilon< \pi/2$ and the tubular neighborhood so that all $\rho>0$ near $r$ in it also satisfy $\rho>1$. A quick computation using \eqref{eq:roots} shows that $t_1, t_2$ stay on one side of the real axis when we move along this tubular neighborhood towards $r$, so we can still separate the compact set $t_{1,2}(T_r)$ from the singularities $\frac{1}{k^2}$, $1 \le k \le d$, as there is no closed loop in $t_{1,2}(T_r)$ encircling one of those. This reduces the inductive determination of the singularities of $F_{d+1}$ to showing that $t_{1,2}(r)$ belongs to the set of singularities of $F_d$ precisely when $r$ is in the claimed set of singularities of $F_{d+1}$ and dealing with the cases $w=1$ and $w=\infty$. The first part is easy since $\frac{r}{(\sqrt r \pm 1)^2}=\frac{1}{k^2}$ if and only if $r=\frac{1}{(k+1)^2}$ or $r=\frac{1}{(k-1)^2}$ (with the usual convention that $1/0=\infty$), so at least outside of the cases $1$ and $\infty$ we indeed get inductively that the singularities of $F_{d+1}$ are included in the set $\frac{1}{k^2}$, $1 \le k \le d+1$, with $k$ and $d+1$ of the same parity.

 A closer analysis of the quadratic polynomial root decomposition, induction, and the analysis in \cite[VI]{FS} show that indeed those are singularities of the same type as the main singularity at $1/(d+1)^2$, so in particular of finite order. This analysis also shows that if we restrict ourselves to the principal branch of $F_d$ defined (for $d \ge 2$) on the plane minus the half line $[\frac{1}{d^2}, \infty]$, we can define inductively $F_{d+1}$ single valued and analytic on $[\frac{1}{(d+1)^2}, \infty]$.

The only issue that remains in the inductive analysis of the singularities of $F_{d+1}$ is what happens when $w=1$ and $w=\infty$, since we have $t_1(\infty)=t_2(\infty)=1$. In the $d\ge 2$ even case, when $1$ is not a singularity of $F_d$, we can first define an analytic function around infinity by the corresponding formula \eqref{eq:ConvExpN} but on a small circle centered at $1$, so
$$
F_{d+1,1}(w)=\frac{1}{2\pi i}\int_{|t-1|=\epsilon} F_d(t)P(t,w)^{-1/2}\, dt, \text{ for } |w|>R_d.
$$

Now there exists a curve $\Gamma$ such that $0$ and $\infty$ belong to $U_{\Gamma}$, but a simple topological argument shows that unfortunately $0$ and $\infty$ will never be in the same connected component of $U_{\Gamma}$. Hence the corresponding integral \eqref{eq:ConvExpNC} will define two separate analytic functions, our $F_{d+1}(w)$ continuing the power series from $0$ on the component of $U_{\Gamma}$ containing $0$ and $F_{d+1,1}(w)$ continuing the power series from $\infty$ on the component of $U_{\Gamma}$ containing $\infty$. Taking curves $\Gamma$ containing zero inside and the singularities of $F_{d+1}$ outside and $\Gamma_1$ containing the singularities of $F_{d+1,1}$ inside, one defines $F_{d+2}$ and $F_{d+2,1}$, respectively, at say any non real positive complex number and then by representing those as double integrals and deforming appropriately, one can show that actually the two integrals define the same function there. Since as above, all points on the extended positive real axis except $1/k^2$, $1 \le k \le d+2, k$ even, $1$ and $\infty$ can be reached by the first construction and all points on the extended positive real axis except $1/k^2$, $1 \le k \le d+2, k$ even, $0$ and $\infty$ can be reached by the second construction, we obtain that $F_{d+2}$ can indeed be extended as a multivalued analytic function to all points except $1/k^2$, $1 \le k \le d+2$, $k$ even, and  $\infty$, and all its singularities are Fuchsian (and of finite order). Note that the analysis before shows the logarithmic nature of the singularity at infinity for all $d \ge 3$. In particular this implies that $1/F$ can never be holonomic for any $d \ge 3$.

Since $F_d$ is holonomic and with singularities of finite order, it also satisfies a Fuchsian ODE as discussed in the previous section, so Theorem~\ref{T:TheoremXs} is proved. We note in Section~\ref{sec:ODEs} that indeed $F_2$, $F_3$, $F_4$, and $F_5$ satisfy ODE's of degree $d-1$ so we conjecture that this holds for all $d \ge 6$.

\vspace{.3cm}
Theorem~\ref{T:TheoremAs} is a simple consequence of the defining relation $A_d(w)=F_d(w)*F_2(w)=F_d(w)*A_1(w)$ and of Theorem~\ref{T:TheoremXs}. This is because the singularities of $F_d$ are finite in number while the singularities of $F_2(w)=A_1(w)=1/\sqrt {1-4w}$ are $1/4$ and $\infty$, so the finite singularities of the convolution are included in the products of those of $F_d$ and $1/4$. An analysis similar to that in \cite[VI]{FS} shows that all of these will be of the same type as the main singularity, since the ones of $F_d$ are so. At infinity, in the odd case we will get only half integer powers as for the finite singularities. In the even case, we get also half integer powers, in addition to the expected logarithmic terms. This happens as in the $A_2$ case and with the same type of analysis at large arguments of the convolution of functions that behave like half integral powers at their finite singularities. $A_1$, $A_2$, $A_3$, $A_4$, and $A_5$ satisfy ODE's of degree $d$ so we conjecture that this holds for all $d \ge 6$ (see Section~\ref{sec:ODEs}).

\vspace{.3cm}
Note that by the theory of Fuchsian ODEs \cite{Ince}, $F_{2d}$ and $A_{2d+1}$ actually have Puiseux expansions in fractional (half integral) powers of $z-w_0$ at all finite singularities, that hold in angular domains near the singularity, while $F_{2d+1}$ and $A_{2d}$ have similar expansions with logarithmic terms instead of fractional powers. So the asymptotic expansions of $\tilde F$ and $\tilde A$ from Theorem~\ref{T:TheoremC} are actually equalities in angular domains at $1$ and extend to similar expansions in the neighborhoods of all the finite singularities and infinity (with the appropriate of an angular domain there).

\vspace{.3cm}
From the relation $(1-B_d)A_d=1$ and the observations that $A_d(1/(2d)^2)=m_d$, with $m_d=\infty$, for $d=1,2$, and $m_d>0$, for $d \ge 3$, and that $B_d$ has positive coefficients and it converges uniformly on the circle of convergence, for all $d \ge 1$ with $|B_d(w)| \le B_d(|w|)=p_d  \le 1$, we conclude that $A_d \ne 0$ on the circle of convergence. All these imply that $B_d=1-1/A_d$ is $\Delta$-analytic too, while the rest of the claims of Theorem~\ref{T:TheoremBs} are clear from the definition of $B$, except for non-holonomicity which will be shown in the next subsection. That is not surprising since $A$ is holonomic but not algebraic for $d \ge 2$ (as will be shown also in the next subsection), so $1/A$ (hence $B$) is then usually non-holonomic.

\section{ODEs and $P$-recursion }\label{sec:ODEs}
We prove that the $x$ and $A$ sequences are generally transcendental over $\mathbb C$, beyond the trivial low dimensional explicit examples, so in particular also $F$ and $A$ are transcendental analytic functions, but they are algebraic (in the sense of generating formal power series discussed earlier) modulo any prime. This is not that surprising considering that ${ {2n \choose n}^2 }$ (which is $A^{(2)}_n$) satisfies this property (\cite{Stanley2}). On the other hand the $B$ sequences are not holonomic, as noted earlier. We also record here the recursions and ODEs for a few low dimensional cases. All the statements above except the non-holonomicity of $B_{2d+1}$, $d \ge 1$, will follow from earlier considerations and known results in the literature. Nevertheless, a careful consideration of the coefficients of the Puiseux expansion of $A'_{2d+1}/A_{2d+1}$ at the first singularity, for example, shows that the first few such terms are in $Q(m_{2d+1}\pi^d)$ so the above methods would require the (unknown but likely) transcendentality of $m_{2d+1}\pi^d$ for the standard proof of the non-algebraicity of $A'_{2d+1}/A_{2d+1}$, and consequently the non-holonomicity of $B_{2d+1}$. Based on the expansion for $A_3$ \cite{Joyce72} and our earlier computations, we conjecture that all the coefficients above are indeed in $Q(m_{2d+1}\pi^d)$. However, both the sequences $( x_n^{(d)} )_n$ and $( A_{2n}^{(d)} )_n$ satisfy the Lucas property \cite{McIntosh}, $u_{np+q}=u_nu_q$ mod $p$, for all $n \ge 0$, $0 \le q \le p-1$ and $p$ prime, which will enable us to give a direct proof that $1/F_d$, $d \ge 3$, and $B_d$, $d \ge 2$, are not holonomic.

\begin{theorem}\label{T:TheoremH}
\begin{enumerate}
\item[(i)] For $d \ge 3$, $( x_n^{(d)} )_n$ is a $P$-recursive sequence, transcendental over $\mathbb C$, but algebraic modulo $p$ prime.
\item[(ii)] For $d \ge 2$, $( A_{2n}^{(d)} )_n$ is a $P$-recursive sequence, transcendental over $\mathbb C$, but algebraic modulo $p$ prime.
\item[(iii)] For $d \ge 2$, $( B_{2n}^{(d)} )_n$ is not a $P$-recursive sequence.
\end{enumerate}
\end{theorem}

We proved in the previous section that the power series $F_d$ and $A_d$ generated by $(x_n^{(d)})_n$ and $(A_{2n}^{(d)})_n $, respectively, are holonomic, which is equivalent with the $P$-recursivity of the given sequences \cite[Ch.6]{Stanley2}. Also the algebraicity modulo any prime $p$ of the sequences follows immediately by induction from the convolution representations of the generating series $F$ and $A$ from the previous section, since the convolution kernels $G$ and $F_2=A_1$ are algebraic (\cite{SW,F}).

Since the singularities of algebraic functions (over $\mathbb C$) cannot contain logarithms, dilogarithms etc. (so the normalized coefficients of the power series at the radius of convergence singularity cannot be asymptotic with negative integral powers of $n$), it also follows that $F_{2d+1}$ and $A_{2d}$, $d \ge 1$, cannot be algebraic over $\mathbb C$ from the results of the previous sections. Similarly, using that holonomic functions cannot have arbitrarily high logarithmic powers or negative logarithmic powers in their asymptotic formulas, it follows from Theorem~\ref{T:TheoremC}~(iii) that $B_{2d}$, $d \ge 1$, cannot be holonomic.

For the case of half odd integral powers asymptotics ($F_{2d}$, $d \ge 2$, and $A_{2d+1}$, $d \ge 1$) transcendentality over $\mathbb C$ follows from arithmetic considerations precisely as in \cite{FGS}, since the integral nature of the coefficients ensure that if one of the corresponding $F_{2d}$, $d \ge 2$, or $A_{2d+1}$, $d \ge 1$, were algebraic over $\mathbb C[x]$, it would actually be algebraic over $\mathbb Q[x]$ by Eisenstein Lemma, hence the Puiseux expansion at the radius of convergence singularity would have algebraic coefficients. However the leading asymptotic coefficient is transcendental as we see in Theorem~\ref{T:TheoremC}, parts (i) and (ii). Note that the coefficients  $\frac{a_2 \, \Gamma(\frac{3-2}{2})}{\pi^{(2-1)/2}}=1$ for $F_2$ and  $\frac{a_1 \, \Gamma(\frac{2-1}{2})}{\pi^{1/2}}=1$ are indeed algebraic as those are the only cases when the $\Gamma$ term cancels the power of $\pi$ in the denominator for either $F$ or $A$ and that of course corresponds to the fact that $F_2=A_1$ is indeed algebraic.

To handle $B_{2d+1}$, $d \ge 1$, we change tack and note that a theorem of Harris and Sibuya \cite{HS} shows that $1/A_{2d+1}$ (hence $B_{2d+1}$) is holonomic if and only if $A'_{2d+1}/A_{2d+1}$ is algebraic. As noted above, an integral sequence $u_n$ is said to have the Lucas property if for all $p$ prime (large enough) $u_{np+q}=u_nu_q$ mod $p$ for all $n \ge 0$, $0 \le q \le p-1$.

Equation \eqref{eq:Xs} and the well known (double) Lucas property \cite{McIntosh} of the binomial coefficients ${np+q \choose kp+r}={n \choose k}{q \choose r}$ mod $p$, where $n,k \ge 0$, $0 \le q,r \le p-1$ and $p$ prime, immediately imply that $x_n^{(d)}$ and $A_{2n}^{(d)}={2n \choose n}  x_n^{(d)}$ have the Lucas property and actually $A_{2n}^{(d)}=0$ mod $p$ when $(p-1)/2 < n \le p-1$. In particular Theorem~6.4 in \cite{Allouche} immediately imply that $F_d$, $d \ge 3$, and $A_d$, $d \ge 2$, are not algebraic over $\mathbb Q$ as they are not of the form $P^{-1/2}$ for a polynomial of degree at most $2$ (from the location and nature of their singularities, for example).

Now let's consider $G(x)=\sum_{n\geq 0} a_n x^n$, where $a_n \in \mathbb Z$, $a_0=1$, satisfy Lucas property and denote $G_p=\sum_{n \le p-1}a_nx^n$. Working modulo $p$, we have
$$
G(x) = \sum_{n \ge 0}\sum_{0 \le q \le p-1} a_{np+q} x^{np} x^q
     = \sum_{n \ge 0}\sum_{0 \le q \le p-1} a_{n} a_q x^{np} x^q
     = G^p(x) G_p(x),
$$
since $G^p(x)=(\sum a_n x^n)^p=\sum a_n^p x^{np}=\sum a_n x^{np}$ by Fermat. Differentiating formally we get
$$
G'(x)=G^p(x)G'_p(x),
$$
since $(G^p)'=0$ mod $p$. Hence $G'/G=G'_p/G_p$ is rational over $\mathbb{Z}/p$.

Next, if $y=G'/G$ is also algebraic over $\mathbb{Q}$ (which is the necessary and sufficient condition for $1/G$ to be holonomic when $G$ is so), there are $k+1$ integral polynomials $Q_0, Q_1, \dots, Q_k$ of degrees at most some fixed $N$ with $Q_k(x) y^k+\dots +Q_0(x)=0$. We can assume the equation is of minimal degree so $Q_0$ and $Q_k$ are not identically zero. Reducing modulo $p$ (for large enough $p$ greater than the absolute value of all the coefficients of all the $Q_m$, so they do not reduce, and $Q_0$ and $Q_k$ are still nonzero mod $p$) and representing $y=G'/G=G'_p/G_p=U_p/V_p$ in lowest terms mod $p$, so $\gcd(U_p,V_p)=1$ in $(\mathbb{Z}/p)[x]$, we get that $U_p|Q_0$ and $V_p|Q_k$ in $(\mathbb{Z}/p)[x]$. This gives that the degrees of $U_p$ and $V_p$ are at most $N$.

Considering now the $2N+2$ integral series $x^m G'$ and $x^m G$, $m=0,1, \dots ,N$, the relation $U_p G-V_p G'=0$ mod $p$ and the fact that the degrees of $U_p$ and $V_p$ are at most $N$ imply that these $2N+2$ integral series are linearly dependent over $\mathbb{Z}/p$, for all large enough primes $p$. Considering now their coefficients as a matrix with $2N+2$ rows and infinitely many columns, we get that any $2N+2$ square sub-matrix has determinant zero mod $p$ for all $p$ large enough, hence is zero as an integer which implies that the $2N+2$ series are linearly dependent over $\mathbb Q$, hence, by clearing denominators, over $\mathbb Z$. So writing the respective linear dependence and grouping the $x^m G'$ together and the $x^m G$ together we get that $P G'+Q G=0$ for some integral polynomials $P$ and $Q$ which are not (both) zero. Hence $G'/G$ must be a rational function over $\mathbb Q$.

To recapitulate, we proved above that a formal power series $G$ with integral coefficients satisfying Lucas property and for which $G'/G$ is algebraic over $\mathbb Q$, must satisfy $G'/G$ be actually rational over $\mathbb Q$. Since clearly none of $F'_d/F_d$, $d \ge 3$, or $A'_d/A_d$, $d \ge 2$, are rational (the first singularity is clearly not a pole, for example), it follows from the Harris-Sibuya theorem (\cite{HS}) that $1/F_d$, $d \ge 3$, and $1/A_d$, $d \ge 2$ (hence $B_d$, $d \ge 2$) cannot be holonomic over $\mathbb C$, or equivalently over $\mathbb Q$ or $\mathbb Z$ (\cite[part VIII, no.151]{PS}).

\vspace{.3cm}
With a more careful analysis, we can actually dispense with the theorem of Harris and Sibuya (or better put, reprove it in our considerably simpler context) and show directly that $1/F_d$ and  $1/A_d$ can be holonomic only in the known low dimensional cases.

In general, if a power series $G$ as above with integral coefficients is holonomic, it is Fuchsian with rational indices (exponents) at all singularities (\cite[6.2]{ABD}). Consequently, if $G$ satisfies Lucas property and $1/G$ is holonomic, our proof showing that $G'/G$ rational forces $G$ to be algebraic implies, using \cite{Allouche}, that $G=P^{-1/2}$ with $P$ an integral polynomial of degree at most $2$. (Indeed, looking at the form $H=\log G$, an integral of a rational function, could have such that $G=e^H$ is Fuchsian with rational indices both at finite singularities and infinity, one concludes that $H$ must be of the form $c+\sum c_k \log (z-z_k)$, with $c$, $c_k \in \mathbb Q$. This implies that $G$ actually satisfies $G^M=R$, for some integer $M$ and rational function $R$.)

\begin{remark}\label{R:x-P-recurrences}
We were able to derive some results that are available in the literature about the sequences $( x_n^{(d)} )_n$. The case $d=3$ is \href{https://oeis.org/A002893}{A002893}, $d=4$ is \href{https://oeis.org/A002895}{A002895}, and $d=5$ is \href{https://oeis.org/A169714}{A169714}, all from \cite{Sloane}, though here all the sequences appear in an unified context (not mentioned in Sloane).

\begin{enumerate}
\item[(i)] The $P$-recurrences for $d=3$, $4$, and $5$, respectively,\footnote{When $d=1$,
$x_{n+1}^{(1)}-x_{n}^{(1)}=0$. When $d=2$, $(n+1) x_{n+1}^{(2)}- 2(2n+1) x_{n}^{(2)}=0$.
}
are:

\begin{equation} \label{eq:P-rec-X-d3}
(n+2)^2 x_{n+2}^{(3)} - (10 n^2 + 30 n + 23) x_{n+1}^{(3)} + 9 (n+1)^2 x_{n}^{(3)}=0;
\end{equation}

\begin{equation} \label{eq:P-rec-X-d4}
(n+2)^3 x_{n+2}^{(4)} - 2 (2n+3) (5 n^2 + 15 n + 12) x_{n+1}^{(4)} + 64 (n+1)^3 x_{n}^{(4)}=0;
\end{equation}

\begingroup\makeatletter\def\f@size{11}\check@mathfonts
\begin{equation} \label{eq:P-rec-X-d5}
 \begin{aligned}
(n & +3)^4 \, x_{n+3}^{(5)} -(35 n^4+350 n^3+1323 n^2 + 2240 n +1433) \, x_{n+2}^{(5)} \\
   & +(n+2)(259 n^3 + 1554 n^2 + 3134 n + 2124) \, x_{n+1}^{(5)}
                  - 3^2\cdot 5^2\cdot (n+1)^2 (n+2)^2 \, x_{n}^{(5)}=0.
 \end{aligned}
\end{equation}
\endgroup

\item[(ii)] The ODEs satisfied by $F_3$, $F_4$ and $F_5$, respectively,
    are given below.\footnote{When $d=1$, the ODE is $(z-1) F' + F=0$. When $d=2$, it is $(4z-1) F' + 2F=0$. } The ODE for the case $n=5$ seems to be new and the $P$ symbols are consistent with the asymptotics obtained in the previous subsections.

\begin{equation} \label{eq:ODE-F-d3}
z(z-1)(9z-1) F''(z)+(27z^2-20z+1) F'(z)+3(3z-1) F(z)=0;
\end{equation}

\begin{equation} \label{eq:ODE-F-d4}
 \begin{aligned}
 z^2 &(4z-1)(16z-1) F'''(z) + 3 z (128z^2-30z+1) F''(z) \\
     &+ (448z^2-68z+1) F'(z)+4(16z-1) F(z)=0.
 \end{aligned}
 \end{equation}

 \begin{equation} \label{eq:ODE-F-d5}
 \begin{aligned}
 z^3 &(z-1)(9z-1)(25z-1) F^{(4)}(z) + z^2 (2700z^3-2590z^2+280z-6) F'''(z)\\
     &+ z (8550z^3-6501z^2+518z-7) F''(z) + (7200z^3-3963z^2+196z+1) F'(z)\\
     &+(900z^2-285z+5) F(z)=0.
 \end{aligned}
 \end{equation}

\end{enumerate}
\end{remark}

\begin{remark}\label{R:A-P-recurrences}
We were also able to derive some results about the sequences
 $( A_n^{(d)} )_n$. The case $d=3$ is \href{https://oeis.org/A002896}{A002896}, $d=4$ is \href{https://oeis.org/A039699}{A039699}, and $d=5$ is \href{https://oeis.org/A287317}{A287317}, all from \cite{Sloane}.
The recurrence and the ODE for the cases $d=4$ and $d=5$ seem to be new; the $P$ symbols are consistent with the asymptotics obtained in the previous subsections.

\begin{enumerate}

\item[(i)] The $P$-recurrences for $d=3$, $4$, and $5$, respectively,\footnote{When $d=1$,
$(n+1) A_{n+1}^{(1)}- 2(2n+1) A_{n}^{(1)}=0$. When $d=2$,
$(n+1)^2 A_{n+1}^{(2)}- 4(2n+1)^2 A_{n}^{(2)}=0$.
} are:

\begingroup\makeatletter\def\f@size{11}\check@mathfonts
\begin{equation} \label{eq:P-rec-A-d3}
(n+2)^3 A_{n+2}^{(3)} - 2 (2n+3)(10 n^2 + 30 n + 23) A_{n+1}^{(3)}
  + 36 (2n+3)(2n+1)(n+1) A_{n}^{(3)}=0;
\end{equation}
\endgroup

\begingroup\makeatletter\def\f@size{11}\check@mathfonts
\begin{equation} \label{eq:P-rec-A-d4}
(n+2)^4 A_{n+2}^{(4)} - 4 (2n+3)^2 (5 n^2 + 15 n + 12) A_{n+1}^{(4)}
  + 256 (2n+3)(2n+1)(n+1)^2 A_{n}^{(4)}=0;
\end{equation}
\endgroup

\begingroup\makeatletter\def\f@size{11}\check@mathfonts
\begin{equation} \label{eq:P-rec-A-d5}
 \begin{aligned}
(n+3)^5 \, A_{n+3}^{(5)} & -2(2n+5)(35 n^4+350 n^3+1323 n^2 + 2240 n +1433) \, A_{n+2}^{(5)} \\
                & +4(2n+5)(2n+3)(259 n^3 + 1554 n^2 + 3134 n + 2124) \, A_{n+1}^{(5)} \\
                & - 1800(2n+5)(2n+3)(2n+1)(n+1)(n+2) \, A_{n}^{(5)}=0.
 \end{aligned}
\end{equation}
\endgroup

\item[(ii)] The ODEs satisfied by $A_3$, $A_4$ and $A_5$, respectively,\footnote{When $d=1$, the ODE is $(4z-1) A' + 2A=0$. When $d=2$, it is $z(16z-1)A''+(32z-1)A'+4A=0$.} are:



\begingroup\makeatletter\def\f@size{11}\check@mathfonts
\begin{equation} \label{eq:ODE-A-d3}
 \begin{aligned}
 z^2 &(4z-1)(36z-1) A'''(z) + 3 z (288z^2-60z+1) A''(z) \\
                            &+ (972z^2-132z+1) A'(z)+6(18z-1) A(z)=0.
 \end{aligned}
\end{equation}
\endgroup

\begingroup\makeatletter\def\f@size{11}\check@mathfonts
\begin{equation} \label{eq:ODE-A-d4}
 \begin{aligned}
 z^3 &(16z-1)(64z-1) A^{(4)}(z) + 2z^2(5120z^2-320z+3)A'''(z) \\
 &+ z (25344z^2-1172z+7) A''(z) + (14592z^2-424z+1) A'(z)+8(96z-1) A(z)=0.
 \end{aligned}
\end{equation}
\endgroup

\begingroup\makeatletter\def\f@size{11}\check@mathfonts
\begin{equation} \label{eq:ODE-A-d5}
 \begin{aligned}
z^4 & (4z-1)(36z-1)(100z-1) A^{(5)}(z) \\
  &+ z^3 (1800*140z^3-62160z^2+1750z-10)A^{(4)}(z)\\
  &+ z^2 (1800*730z^3-268740z^2+5992z-25)A'''(z) \\
  &+ z (1800*1275z^3-369240z^2+5964z-15) A''(z) \\
  &+ (1800*600z^3-124020z^2+1196z-1) A'(z) \\
  &+(1800*30z^2-3420z+10) A(z)=0.
 \end{aligned}
\end{equation}
\endgroup

\end{enumerate}
\end{remark}


\end{document}